\documentclass[11pt,a4paper,twoside]{article}
\usepackage{graphicx}
\usepackage{amsmath}
\usepackage{amsfonts}
\usepackage{amssymb}
\usepackage{enumerate}
\usepackage[hidelinks]{hyperref}
\usepackage{lscape}
\usepackage{longtable}
\usepackage{mathrsfs}
\usepackage{rotating}
\usepackage{color}
\usepackage{url}
\usepackage{rotating}
\usepackage{yfonts}

\usepackage[labelfont=sl,textfont=sl]{subcaption} % newer version of subfig
\usepackage[font={small,sl},labelsep=period]{caption} % format figure & table captions

\usepackage{titlesec} 
\usepackage[titletoc]{appendix}
\usepackage[ruled,vlined,noline,linesnumbered]{algorithm2e}
\newtheorem{theorem}{Theorem}[section]

\newtheorem{corollary}{Corollary}[section]

\newtheorem{definition}{Definition}[section]
\newtheorem{example}{Example}[section]

\newtheorem{lemma}{Lemma}[section]

\newtheorem{proposition}{Proposition}[section]
\newtheorem{remark}{Remark}[section]

\newenvironment{proof}[1][Proof]{\textbf{#1.} }{\hfill$\Box$}
\DeclareMathOperator{\diag}{diag}
\DeclareMathOperator{\spann}{span}
\DeclareMathOperator{\pspan}{pspan}

\newcommand{\T}{\top}
\DeclareMathOperator{\argminop}{arg\,min}
\newcommand{\argmin}[1]{\underset{#1}{\argminop}}
\DeclareMathOperator{\argmaxop}{arg\,max}
\newcommand{\argmax}[1]{\underset{#1}{\argmaxop}}
\newcommand{\R}{\mathbb{R}}
\newcommand{\LL}{\mathbb{L}}
\newcommand{\cm}[1]{\operatorname{cm}\left(#1\right)}
\newcommand{\cmk}[2]{\operatorname{cm}_{#2}\left(#1\right)}
\newcommand{\D}{\mathcal{D}}
\newcommand{\Dm}{\mathbf{D}}

\newcommand{\cV}{c\mathcal{V}}

\newcommand{\bbm}{\begin{bmatrix}}
\newcommand{\ebm}{\end{bmatrix}}

\newcommand{\vd}{\mathbf{d}}
\newcommand{\vc}{\mathbf{c}}
\newcommand{\ve}{\mathbf{e}}
\newcommand{\vu}{\mathbf{u}}
\newcommand{\vv}{\mathbf{v}}

\newcommand{\zero}{\mathbf{0}}
\newcommand{\zeron}{\mathbf{0}_n}
\newcommand{\one}{\mathbf{1}}
\newcommand{\onen}{\mathbf{1}_n}

\newcommand{\G}{\mathbf{G}} % Gram matrix

\newcommand{\B}{\mathcal{B}}
\newcommand{\setI}{\mathcal{I}}
\newcommand{\In}{\mathcal{I}_n}
\newcommand{\mIn}{\mathbf{I}_n}
\newcommand{\mI}{\mathbf{I}}

\newcommand{\Pm}{\mathbf{P}}
\newcommand{\Rm}{\mathbf{R}}
\newcommand{\Sm}{\mathbf{S}}
\newcommand{\Bm}{\mathbf{B}}

\newcommand{\set}{\mathcal{S}}
\newcommand{\setS}{\mathcal{S}}

\newcommand{\dlk}{\mathcal{D}_{\LL}^{(1:k)}}

\def\sec52{0} % Can remove everything inbetween the if conditions

\begin{document}

\title{\centering Using orthogonally structured positive bases for constructing
positive $k$-spanning sets with cosine measure guarantees}
\author{
Warren Hare 
\thanks{Department of Mathematics, University of British Columbia, Kelowna, 
British Columbia, Canada.  Hare's research is partially funded by the 
Natural Sciences and Engineering Research Council of Canada (cette recherche est partiellement financ\'ee par le 
Conseil de recherches en sciences naturelles et en g\'enie du Canada),
Discover Grant \#2018-03865. ORCID 0000-0002-4240-3903 
(\texttt{warren.hare@ubc.ca}).}           %  \\
\and
Gabriel Jarry-Bolduc
\thanks{Mathematics and Statistics Department, Saint Francis Xavier Univeristy, Antigonish, 
Nova Scotia, Canada.  Jarry-Bolduc's research is partially funded 
through the Natural Sciences and Engineering Research Council of 
Canada  (cette recherche est partiellement financ\'ee par le 
Conseil de recherches en sciences naturelles et en g\'enie du Canada), Discover Grant \#2018-03865. ORCID 0000-0002-1827-8508 
(\texttt{gabjarry@alumni.ubc.ca}).}
\and
S\'ebastien Kerleau\thanks{LAMSADE, CNRS, Universit\'e Paris Dauphine-PSL, 
Place du Mar\'echal de Lattre de Tassigny, 75016 Paris, France 
(\texttt{sebastien.kerleau@lamsade.dauphine.fr}).} 
\and 
Cl\'ement W. Royer\thanks{LAMSADE, CNRS, Universit\'e Paris Dauphine-PSL, 
Place du Mar\'echal de Lattre de Tassigny, 75016 Paris, France. Royer's 
research is partially funded by Agence Nationale de la Recherche through 
program ANR-19-P3IA-0001 (PRAIRIE 3IA Institute). ORCID 0000-0003-2452-2172
(\texttt{clement.royer@lamsade.dauphine.fr}).}
\thanks{
Research for this paper is also supported by France-Canada Research Funds 
2022.
}
}

\date{\today}

%\titlerunning{OSPBs for uilding P$k$SS}

\maketitle

\begin{abstract}
Positive spanning sets span a given vector space by nonnegative linear 
combinations of their elements. These have attracted significant attention in 
recent years, owing to their extensive use in derivative-free optimization. In 
this setting, the quality of a positive spanning set is assessed through its 
cosine measure, a geometric quantity that expresses how well such a set covers 
the space of interest. 
In this paper, we investigate the construction of positive $k$-spanning sets 
with geometrical guarantees.
Our results build on recently identified positive spanning sets, called orthogonally structured positive bases.
We 
first describe how to identify such sets and compute their cosine measures 
efficiently. We then focus our study on positive $k$-spanning sets, for which 
we provide a complete description, as well as a new notion of cosine measure 
that accounts for the resilient nature of such sets. By combining our 
results, we are able to use orthogonally structured positive bases to 
create positive $k$-spanning sets with guarantees on the value of their cosine 
measures.
\end{abstract}

 
\paragraph{2020 Mathematics Subject Classification} 
% 15A03 Vector spaces, linear dependence, rank, lineability
% 15A21 Canonical forms, reduction, classification
% 15B30 Orthogonal matrices
% 15B99 Special matrices (none of the specific categories)
% 90C56 Derivative-free methods and methods using generalized derivatives
15A03; 15A21; 15B30; 15B99; 90C56.

%%%%%%%%%%%%%%%%%%%%%%%%%%%%%%%%%%%%%%%%%%%%%%%%%%%%%%%%%%%%%%%%%%%%%%%%%%%%%%%
\section{Introduction}
\label{sec:intro}
%%%%%%%%%%%%%%%%%%%%%%%%%%%%%%%%%%%%%%%%%%%%%%%%%%%%%%%%%%%%%%%%%%%%%%%%%%%%%%%

Positive spanning sets, or PSSs, are sets of vectors that span a space 
of interest through nonnegative linear combinations. These sets were first 
investigated by Davis~\cite{Davis1954positive}, with a particular emphasis on 
inclusion-wise minimal PSSs, also called positive 
bases~\cite{Davis1954positive,mckinney1962}. Such a restriction guarantees 
bounds on the minimal and maximal size of a positive 
basis~\cite{Audet2011short,Davis1954positive}. PSSs of this form are now well 
understood in the literature~\cite{Regis2016}. Any positive basis is amenable 
to a decomposition into minimal positive bases over 
subspaces~\cite{Romanowicz1987}, a result that recently lead to the 
characterization of nicely structured positive bases~\cite{hare2023nicely}, 
hereafter called orthogonally structured positive bases, or OSPBs.

Positive spanning sets and positive bases are now a standard tool to develop 
derivative-free optimization algorithms~\cite{Hare2018,Conn2009}. The idea is 
to use directions from a PSS in replacement for derivative information. In 
that setting, it becomes critical to quantify how well directions from a PSS 
cover the space of interest, which is typically assessed through the cosine 
measure of this PSS~\cite{lewis1996rank,kolda2003directsearch}. In general, no 
closed form expression is available for this cosine measure, yet various 
computing techniques have been proposed~\cite{hare2020cm,Regis2021}. In 
particular, a deterministic algorithm to compute the cosine measure of any PSS 
was recently proposed by Hare and Jarry-Bolduc~\cite{hare2020cm}. The cosine 
measure of certain positive bases is known, along with upper bounds for 
positive bases of minimal and maximal sizes~\cite{naevdal2019,Regis2016}. The 
intermediate size case is less understood, even though recent progress was 
made in the case of orthogonally structured positive 
bases~\cite{hare2023nicely}. 

Meanwhile, another category of positive spanning sets introduced in the 
1970s has received little attention in the literature. Those sets, called 
positive $k$-spanning sets (P$k$SSs), can be viewed as resilient PSSs, since 
at least $k$ of their elements must be removed in order to make them lose 
their positively spanning property~\cite{Marcus1981}. Unlike standard PSSs, 
only partial results are known regarding inclusion-wise minimal P$k$SS. These 
results depart from those for PSSs as they rely on polytope 
theory~\cite{Marcus1984,Wotzlaw2009}. Moreover, the construction of P$k$SSs 
based on PSSs has not been fully explored. To the best of our knowledge, the 
cosine measure of P$k$SSs has not been investigated, which partly prevents 
those sets from being used in derivative-free algorithms.

% Contributions
In this paper, we investigate positive $k$-spanning sets through the lens of 
orthogonally structured positive bases and cosine measure. To this end, we 
refine previous results on OSPBs so as to obtain an efficient cosine measure 
calculation technique. We then explain how the notion of cosine measure can 
be generalized to account for the positive $k$-spanning property, thereby 
introducing a quantity called the $k$-cosine measure. To the best of our 
knowledge, this definition is new in both the derivative-free optimization 
and the positive spanning set literature. By combining those elements, we are 
able to build positive $k$-spanning sets from OSPBs as well as to provide a 
bound on their $k$-cosine measures. Our results pave the 
way for using positive $k$-spanning sets in derivative-free algorithms.

% Outline
The rest of this paper is organized as follows. Section~\ref{sec:defs}
summarizes important properties of positive spanning sets, including their 
characterization through the cosine measure, with a particular 
focus on the subclass of orthogonally structured positive bases (OSPBs).
Section~\ref{sec:findOSPB} describes an efficient way to compute the cosine 
measure of an OSPB based on leveraging its decomposition. 
Section~\ref{sec:pksscmk} formalizes the main properties associated with 
positive $k$-spanning sets, introduces the $k$-cosine measure to help 
studying these sets and uses OSPBs to design P$k$SSs 
with guaranteed $k$-cosine measure. Section~\ref{sec:conc} summarizes our 
findings and provides several perspectives of our work.

%%%%%%%%%%%%%%%%%%%%%%%%%%%%%%%%%%%%%%%%%%%%%%%%%%%%%%%%%%%%%%%%%%%%%%%%%%%%%%%

\paragraph{Notations}
Throughout this paper, we will work in the Euclidean space $\R^n$ with 
$n \ge 2$, or a linear subspace thereof, denoted by $\LL \subset \R^n$. 
More generally, blackboard bold letters will be used to designate 
infinite-valued sets, such as linear (sub)spaces and half-spaces. 
The Minkowski sum $\left\{\mathbf{a}+\mathbf{b},
(\mathbf{a},\mathbf{b})\in\mathbb{A}\times\mathbb{B}\right\}$ of two spaces 
$\mathbb{A}$ and $\mathbb{B}$ will be denoted by $\mathbb{A}+\mathbb{B}$ or by 
$\mathbb{A}\perp \mathbb{B}$ if the two spaces are orthogonal to one another. 
The orthogonal to a given subspace $\mathbb{A}$ will be denoted accordingly by 
$\mathbb{A}^{\perp}$.
In order to allow for repetition among their elements, positive spanning sets 
will be seen as families of vectors. We will use calligraphic letters such as 
$\D$ to denote finite families of vectors. Given a family $\D \subset \R^n$, 
the linear span of this family ($i.e.$ the set of linear combinations of its 
elements) will be denoted by $\spann(\D)$. Bold lowercase (resp. uppercase) 
letters will be used to designate vectors (resp. matrices). The notations 
$\zeron$ and $\onen$ will respectively be used to designate the null vector 
and the all-ones vector in $\R^n$, while $\In :=\{\ve_1,\cdots,\ve_n\}$ will 
denote the family formed by coordinate basis vectors in $\R^n$. Given a family 
$\D_a$ of vectors in $\R^n$, we will use $\Dm_a$ to denote a matrix with $n$ 
rows and whose columns correspond to the elements in $\D_a$, in no particular 
order unless otherwise stated. We will use $\D_a\setminus\{\vd\}$ to denote a family obtained from $\D_a$ by removing exactly \textbf{one} instance of the vector $\vd$. Similarly, the notation $\D_a\cup\{\vd\}$ will designate the family obtained from  $\D_a$ by adding one copy of the vector $\vd .$ For instance, $\{\vd,\vd,\vd'\}\setminus\{\vd\}=\{\vd,\vd'\}$ and $\{\vd,\vd'\}\cup\{\vd\}=\{\vd,\vd,\vd'\}$.   Finally, the notation $[\![1,m]\!]$ will refer to 
the set $\{1,\dots,m\}$.

%%%%%%%%%%%%%%%%%%%%%%%%%%%%%%%%%%%%%%%%%%%%%%%%%%%%%%%%%%%%%%%%%%%%%%%%%%%%%%%
\section{Background on positive spanning sets}
\label{sec:defs}
%%%%%%%%%%%%%%%%%%%%%%%%%%%%%%%%%%%%%%%%%%%%%%%%%%%%%%%%%%%%%%%%%%%%%%%%%%%%%%%

In this section, we introduce the main concepts and results on positive 
spanning sets that will be used throughout the paper. Section~\ref{subsec:pss} 
defines positive spanning sets and positive bases. 
Section~\ref{subsec:ospbdef} introduces the concept of orthogonally structured 
positive bases, that were first studied under a different 
name~\cite{hare2023nicely}. Section~\ref{subsec:cm} provides the definition of 
the cosine measure of a positive spanning set, along with its value for 
several commonly used positive bases.

%%%%%%%%%%%%%%%%%%%%%%%%%%%%%%%%%%%%%%%%%%%%%%%%%%%%%%%%%%%%%%%%%%%%%%%%%%%%%%%
\subsection{Positive spanning sets and positive bases}
\label{subsec:pss}

In this section, we recall the definitions of positive spanning sets and 
positive bases, based on the seminal paper of 
Davis~\cite{Davis1954positive}. 

\begin{definition}[Positive span and positive spanning set]\label{def:pspanpss}
	Let $\LL$ be a linear subspace of $\R^n$ and $m \ge 1$. The 
	\emph{positive span} of a family of vectors 
	$\D=\left\{\vd_1, \vd_2,\dots, \vd_m \right\}$ in $\LL$, 
	denoted $\pspan(\D)$, is the set 
	$$
		\pspan(\D):=\{\vu \in \LL: 
		\vu=\alpha_1\vd_1+\cdots+\alpha_m \vd_m \text{ }| \text{ }
		\forall i \in [\![1, m]\!], \alpha_i \geq 0\}.
	$$
	A \emph{positive spanning set} (PSS) of $\LL$  is a family of vectors $\D$ 
	such that $\pspan(\D)=\LL$.
\end{definition}
When $\LL$ is not specified, a positive spanning set is understood as 
positively spanning $\R^n$. Throughout the paper, we only consider positive 
spanning sets formed of nonzero vectors. Note, however, that we do not force 
all elements of a positive spanning set to be distinct.
The next two lemmas are well known and describe basic properties of positive spanning sets.

\begin{lemma}\textup{\cite[Theorem 2.5]{Regis2016}} \label{lem:proppss}
Let $\D$ be a finite set of vectors in a subspace $\LL$ of $\R^n$. The following statements are equivalent:
\begin{enumerate}[(i)]
    \item $\D$ is a PSS of $\LL$.
    \item\label{lem: PSS eq pos dotprod}For every nonzero vector 
    $\vu \in \mathbb{L}$, there exists an element $\vd \in \D$ such that 
    $\vu^\T \vd>0$.
    \item $span(\D)=\LL$ and $\zeron$ 
    can be written as a \textbf{positive} linear combination of the elements 
    of $\D$. 
    \end{enumerate}
\end{lemma}
\begin{lemma}\textup{\cite[Theorem 3.7]{Davis1954positive}}\label{lem:pss => linspanning}
    If $\D$ is a PSS of $\LL$, then for any $\vd \in \D$ the set 
    $\D\backslash\{\vd\}$ linearly spans $\LL$.
\end{lemma}

Note that Lemma $\ref{lem:pss => linspanning}$ implies that a PSS 
must contain at least $\dim(\LL)+1$ vectors. There is no upper bound 
on the number of elements that a PSS may contain, but such a restriction 
holds for a subclass of positive spanning sets called positive bases.

\begin{definition}[Positive basis] \label{def:pb}
Let $\LL$ be a linear subspace of $\R^n$ of dimension $\ell \ge 1$.
	A \emph{positive basis} of $\LL$ of size $\ell+s$, denoted $\D_{\LL,s}$, is 
	 a PSS of $\LL$ with $\ell+s$ elements satisfying:
  
  $$\forall \mathbf{d} \in \D_{\LL,s}, \quad\pspan{(\D_{\LL,s}\backslash\{\mathbf{d}\})} \neq \LL.$$ When 
	$\LL=\R^n$, such a set will be denoted by $\D_{n,s}$.
\end{definition}
In other words, positive bases of $\LL$ are inclusion-wise minimal positive spanning sets of $\LL$~\cite{Conn2009}. Positive bases can also be defined thanks to the notion of positive independence.
\begin{definition}[Positive independence]
\label{def:pi}
	Let $\LL$ be a linear subspace of $\R^n$ of dimension $\ell \ge 1$. A family 
	of vectors $\D$ in $\LL$ is \emph{positively independent} if, for any 
	$\vd \in \D$, there exists a vector $\vu \in \LL$ such that 
	$\vu^\T \vd > 0$ and $\vu^\T \vv \le 0$ for any 
	$\vv \in \D \setminus \{\vd\}$.
\end{definition}
A positive basis of $\LL$ is thus a positively independent PSS of $\LL$.In the 
case $\LL=\R^n$, we simply say that $\D_{n,s}$ is a positive basis (of size 
$n+s$).

One can show that the size of a positive basis of a subspace $\LL$ 
is at least $\dim(\LL)+1$ and at most 
$2\dim(\LL)$~\cite{AudetHare2018derivative,Davis1954positive}.

\begin{definition}\label{def:pbminmax}
	Let $\LL$ be a linear subspace of $\R^n$ of dimension $\ell \ge 1$ and 
	$1 \le s \le \ell$. A positive basis $\D_{\LL,s}$ of $\LL$ is called 
	\emph{minimal} if $s=1$, in which case we denote it by 
	$\D_{\LL}=\D_{\LL,1}$. The positive basis is called \emph{maximal} if 
	$s=\ell$. If $1 < s < \ell$, we say that the positive basis has 
	\emph{intermediate size}.
\end{definition}

Maximal and minimal positive bases have a well-understood 
structure~\cite{AudetHare2018derivative,Regis2015}. The structure of an 
arbitrary positive basis can however be quite complicated to analyze. In the 
next section, we describe a decomposition formula for positive bases that 
identifies a favorable structure.

%%%%%%%%%%%%%%%%%%%%%%%%%%%%%%%%%%%%%%%%%%%%%%%%%%%%%%%%%%%%%%%%%%%%%%%%%%%%%%%%
\subsection{Orthogonally structured positive bases}
\label{subsec:ospbdef}

A complete characterization of positive bases can be provided through 
a subspace decomposition~\cite{Romanowicz1987}. Such a decomposition is based 
upon the concept of critical vectors defined below.

\begin{definition}[Critical vectors] 
\label{def:criticalset}
	Let $\LL$ be a subspace of $\R^n$ and $\D_{\LL,s}$ be a positive basis of 
	$\LL$. A vector $\vc \in \LL$ is called a \emph{critical vector} of 
	$\D_{\LL,s}$ if it cannot replace an element of $\D_{\LL,s}$ to form a 
	positive spanning set, i.e.
	$$
	%\label{eq:criticalvec}
		\forall \vd \in \D_{\LL,s},\quad 
		\pspan\left((\D_{\LL,s} \setminus \{\vd\}) \cup \{\vc\} \right) 
		\neq \LL.
	$$
\end{definition}

Note that the zero vector is a critical vector for every positive basis of 
$\LL$, and that it is the only critical vector for a maximal positive 
basis~\cite{Romanowicz1987}. Moreover, if $\D_{\LL}$ is a minimal positive 
basis and $\dim(\LL) \ge 2$, the set of critical vectors (known as the 
complete critical set) is given by
\[
	-\bigcup_{i \neq j} \pspan\left( \D_{\LL} \setminus {\{\vd_i,\vd_j\}} \right).
\]
Using critical vectors, one may decompose any positive basis as a union 
of minimal positive bases. 
Theorem~\ref{thm:decomppb} gives the result for a positive basis in $\R^n$, 
by adapting the original decomposition result due to 
Romanowicz~\cite[Theorem 1]{Romanowicz1987} (see 
also~\cite[Lemma 25]{hare2023nicely}).

\begin{theorem}[Structure of a positive basis] \label{thm:decomppb}
	Suppose that $n \ge 2$, $s\ge 1$, and consider a positive basis $\D_{n,s}$ 
	of $\R^n$. Then, either $s=1$ or there exist subspaces $\LL_1,\dots,\LL_s$ of $\R^n$ 
	such that $\R^n=\LL_1 + \LL_2 + \cdots + \LL_s$ and 
	$\LL_i \cap \LL_j = \{\zero\}$ for any $i\neq j$, and there exist $s$ 
	associated minimal positive bases $\D_{\LL_1},\dots,\D_{\LL_s}$ such that  
	\begin{equation}
	\label{eq:decomppb}
    	\D_{n,s} \; = \; 
    	\D_{\LL_1} \cup (\D_{\LL_2} + \vc_1) \cup 
    	\cdots \cup (\D_{\LL_s} + \vc_{s-1}),
	\end{equation}
	where for any $j=1,\dots,s-1$, we let 
	$\D_{\LL_{j+1}} + \vc_j := \left\{ \vd + \vc_j\ \middle|\ 
	\vd \in \D_{\LL_{j+1}} \right\}$, and the vector 
	$\vc_j \in \LL_1 + \cdots + \LL_j$ 
	is a critical vector for $\D_{\LL_1}\cup\bigcup\limits_{i=2}^{j} (\D_{\LL_i}+ \vc_{i-1})$.
\end{theorem}

The result of Theorem~\ref{thm:decomppb} is actually an equivalence, in that 
any set that admits the decomposition~\eqref{eq:decomppb} is a non-minimal 
positive basis of $\R^n$~\cite[Theorem 1]{Romanowicz1987}. The example below 
provides an illustration for both Definition~\ref{def:criticalset} and 
Theorem~\ref{thm:decomppb}. One can easily check that the stated vector is 
indeed critical, and that the proposed decomposition 
matches~\eqref{eq:decomppb}.

\begin{example} \label{ex:exdecomppb}
	Consider the following positive basis of $\R^4$:
	\[
		\D_{4,2}=\left\{ \bbm 1 \\ 0 \\0 \\0 \ebm,
		\bbm 0 \\ 1 \\ 0 \\0 \ebm,
		\bbm -1 \\ -1\\ 2 \\ 2 \ebm,
		\bbm 1 \\1 \\ -4 \\ -4\ebm 
		\bbm 0 \\ 0 \\ 1 \\ 0 \ebm,
		\bbm 0 \\ 0 \\ 0 \\ 1 \ebm
		\right\}.
	\]
	The first four vectors of this set form a minimal positive basis for 
	$\{[x\ y\ z\ z]^\T , (x,y,z)\in \R^3\}$ that admits 
	$[0\ 0\ 0.5\ 0.5]^\T$ as a critical vector. Therefore, a decomposition of 
	the form~\eqref{eq:decomppb} for $\D_{4,2}$ is:
	\[
		\D_{4,2} = \left\{
		\bbm 1 \\ 0 \\0 \\0 \ebm,
		\bbm 0 \\ 1 \\ 0 \\0 \ebm,
		\bbm -1 \\ -1\\ 2 \\ 2 \ebm,
		\bbm 1 \\1 \\ -4 \\ -4\ebm \right\} 
		\bigcup
		\left(
		\left\{
		\bbm 0 \\ 0 \\ 0.5 \\ -0.5 \ebm,
		\bbm 0 \\ 0 \\ -0.5 \\ 0.5 \ebm
		\right\}
		+ \left\{\bbm 0 \\ 0 \\ 0.5 \\ 0.5 \ebm \right\}
		\right).
	\]
\end{example}

In general, however, the decomposition~\eqref{eq:decomppb} is hard to compute, 
as it requires to determine the critical vectors and the subspaces $\LL_i$, 
that need not be orthogonal to one another. These considerations have lead 
researchers to consider a subclass of positive bases with a nicer 
decomposition~\cite{hare2023nicely}, leading to Definition~\ref{def:ospb} 
below.

\begin{definition}[Orthogonally structured positive bases]
\label{def:ospb}
	Suppose that $n \ge 2$, and consider a positive basis $\D_{n,s}$ of $\R^n$. 
	If there exist subspaces $\LL_1,\dots,\LL_s$ of $\R^n$ such that 
	$\R^n=\LL_1 \perp \LL_2 \perp \cdots \perp \LL_s$ and associated 
	minimal positive bases $\D_{\LL_1},\dots,\D_{\LL_s}$ such that  
	\begin{equation}
	\label{eq:decompospb}
    	\D_{n,s} \; = \; 
    	\D_{\LL_1} \cup \D_{\LL_2} \cup \cdots \cup \D_{\LL_s},
	\end{equation}
	then $\D_{n,s}$ is called an \emph{orthogonally structured positive basis}, 
	or OSPB.
\end{definition}

By construction, note that for any $1 \le i<j \le s$, any pair of elements 
$(\vd,\vd') \in \D_{\LL_i} \times \D_{\LL_j}$ satisfies $\vd^\T \vd'=0$.

The class of orthogonally structured positive bases includes common positive 
bases, such as minimal ones and certain maximal positive bases such as those 
formed by the coordinate vectors and their negatives in $\R^n$. More OSPBs
can be obtained via numerical procedures, even for intermediate 
sizes~\cite{hare2023nicely}. As will be shown in Section~\ref{sec:findOSPB}, 
it is also possible to compute their cosine measure in an efficient manner.

%%%%%%%%%%%%%%%%%%%%%%%%%%%%%%%%%%%%%%%%%%%%%%%%%%%%%%%%%%%%%%%%%%%%%%%%%%%%%%%
\subsection{The cosine measure}
\label{subsec:cm}

To end this background section, we define the cosine measure, a metric 
associated with a given family of vectors.

\begin{definition}[Cosine measure and cosine vector set]
\label{def:cm}
	Let $\D$ be a nonempty family of nonzero vectors in $\R^n$. The \emph{cosine measure} 
	of $\D$ is given by
	$$
		\cm{\D} = \min_{\substack{\|\vu\| =1 \\ \vu \in \R^n}} 
		\max_{\vd \in \D} \frac{\vu^\T \vd}{\|\vd\|}.
	$$
	The \emph{cosine vector set} associated with $\D$ is given by
	\[
		\cV(\D) = \argmin{\substack{\|\vu\| =1 \\ \vu \in \R^n}} 
		\max_{\vd \in \D} \frac{\vu^\T \vd}{\|\vd\|}.
	\]
\end{definition}

The cosine measure and cosine vector set provide insights on the 
geometry of the elements of $\D$ in the space. Such concepts are of particular 
interest in the case of positive spanning sets, as shown by the result of 
Proposition~\ref{prop:cmpss}. Its proof can be found in key 
references in derivative-free optimization~\cite{AudetHare2018derivative,
Conn2009}, but note that our analysis in Section~\ref{sec:pksscmk} provides a 
more general proof in the context of positive $k$-spanning sets.

\begin{proposition}
\label{prop:cmpss}
	Let $\D$ be a nonempty, finite family of vectors in $\R^n$. Then $\D$ is a 
	positive spanning set in $\R^n$ if and only if $\cm{\D}{}>0$.
\end{proposition}

Proposition~\ref{prop:cmpss} shows that the positive spanning property can be 
checked by computing the cosine measure. The value of the cosine measure 
further characterizes how well that set covers the space, which is a relevant 
information in derivative-free algorithms~\cite{kolda2003directsearch}. Both 
observations motivate the search for efficient techniques to compute the 
cosine measure.

We end this section by providing several examples of cosine measures for 
orthogonally structured positive bases, that form the core interest of this 
paper. We focus on building those positive bases from the coordinate vectors. 
A classical choice for a maximal positive basis consists in selecting 
the coordinate vectors and their negatives. In that case, it is known~\cite{naevdal2019} that
\[
	\cm{\In \cup -\In}
	=\cm{\left\{\ve_1,\dots,\ve_n,-\ve_1,\dots,-\ve_n\right\}}{}
	=\frac{1}{\sqrt{n}}.
\]
One can also consider a minimal positive basis formed by the coordinate 
vectors and the sum of their negatives. Then, we have
\[
	\cm{\In \cup \left\{-\sum_{i=1}^n \ve_i\right\}}
	= \cm{\left\{\ve_1,\dots,\ve_n,-\sum_{i=1}^n \ve_i\right\}}{} 
	= \frac{1}{\sqrt{n^2+2(n-1)\sqrt{n}}}.
\]
To the best of our knowledge, the latter formula has only 
been recently established in \cite{jarrybolducthesis,roberts2022subspace}, and 
no formal proof is available in the literature. In the next section, we will 
develop an efficient cosine measure calculation technique for orthogonally 
structured positive bases that will provide a proof of this result as a 
byproduct.

%%%%%%%%%%%%%%%%%%%%%%%%%%%%%%%%%%%%%%%%%%%%%%%%%%%%%%%%%%%%%%%%%%%%%%%%%%%%%%%
\section{Computing OSPBs and their cosine measure}
\label{sec:findOSPB}
%%%%%%%%%%%%%%%%%%%%%%%%%%%%%%%%%%%%%%%%%%%%%%%%%%%%%%%%%%%%%%%%%%%%%%%%%%%%%%%

In this section, we describe how orthogonally structured positive bases can be 
identified using the decomposition introduced in the previous section. We also 
leverage this decomposition to design algorithms for detecting the OSPB 
structure and computing the cosine measure of a given OSPB.

%%%%%%%%%%%%%%%%%%%%%%%%%%%%%%%%%%%%%%%%%%%%%%%%%%%%%%%%%%%%%%%%%%%%%%%%%%%%%%%
\subsection{Structure and detection of OSPBs}
\label{subsec:structOSPB}

The favorable structure of an OSPB can be revealed through the properties of 
its matrix representations, and in particular the Gram matrices associated 
with an OSPB, in the sense of Definition~\ref{de:grammat} below.

\begin{definition}[Gram matrix]
\label{de:grammat}
	Let $\set$ be a finite family of $m$ vectors in $\R^n$ and 
	$\Sm \in \R^{n \times m}$ a matrix representation of this family. The 
	Gram matrix of $\set$ associated with $\Sm$ is 
	the matrix $\G(\Sm)=\Sm^\T \Sm$.
\end{definition}

Given a Gram matrix $\G(\Sm)$ of $\set$, any Gram matrix of $\set$ has the 
form $\Pm^\T \G(\Sm)\Pm$ where $\Pm$ is a permutation matrix. We will 
show that when $\set$ is an OSPB, there exists a Gram matrix with a 
block-diagonal structure that reveals the decomposition of the basis, and 
thus its OSPB nature.

\begin{theorem} 
\label{thm:ospbblocks}
	Let $n \ge 2$, $s\ge 2$ and $\D_{n,s} \subset \R^n$ be a positive basis. Then, 
	$\D_{n,s}$ is orthogonally structured if and only if one of its Gram matrices 
	$\G(\D_{n,s})$ can be written as a block diagonal matrix with $s$ diagonal 
	blocks.
\end{theorem}
\begin{proof}
	Suppose first that $\D_{n,s}$ is an OSPB associated to the decomposition
	$$
		\D_{n,s}=\D_{\LL_1} \cup \D_{\LL_2} \cup \cdots \cup \D_{\LL_{s}}.
	$$
	Let $\Dm_{n,s}$ be a matrix representation of $\D_{n,s}$ such that 
	$\Dm_{n,s}=\bbm \Dm_{\LL_1} &\cdots &\Dm_{\LL_s} \ebm$, where $\Dm_{\LL_i}$ 
	is a matrix representation of $\D_{\LL_i}$ for every $i \in [\![1,s]\!]$. By 
	orthogonality of those matrices, the Gram matrix 
	$\G(\Dm_{n,s})=\Dm_{n,s}^\T \Dm_{n,s}$ is then block diagonal with $s$ 
	blocks corresponding to $\G(\Dm_{\LL_1}),\dots,\G(\Dm_{\LL_s})$, and thus 
	the desired result holds.
	
	Conversely, suppose that there exists a matrix representation $\Dm_{n,s}$ 
	of $\D_{n,s}$ such that $\G(\Dm_{n,s})$ is block diagonal with $s$ 
	blocks. Then, by considering a partitioning with respect to those diagonal 
	blocks, one can decompose $\Dm_{n,s}$ into $\bbm \Sm_1 &\cdots &\Sm_s \ebm$ 
	so that $\G(\Dm_{n,s})=\diag\left( \G(\Sm_1),\dots, \G(\Sm_{s})\right)$. By 
	definition of the Gram matrix, this structure implies that the columns of 
	$\Sm_i$ and $\Sm_j$ are orthogonal for any $i \neq j$, thus we only need 
	to show that each $\Sm_i$ is a minimal positive basis for its linear 
	span. 
	To this end, we use the fact that $\D_{n,s}$ is positively spanning. By 
	Lemma~\ref{lem:proppss}(ii), there exists a vector $\vu \in \R^{n+s}$ with 
	positive coefficients such that $\Dm_{n,s}\vu=\zeron$. By decomposing the 
	vector $\vu$ into $s$ vectors $\vu_1,\dots,\vu_s$ according to the 
	decomposition $\bbm \Sm_1 &\cdots &\Sm_s \ebm$, we obtain that
	\begin{equation}
	\label{eq:ospbth}
		\Dm_{n,s}\vu = \sum_{i=1}^n \Sm_i \vu_i = \zeron.
	\end{equation}
	By orthogonality of the columns of the $\Sm_i$ matrices, the property~\eqref{eq:ospbth} 
	is equivalent to
	\[
		\Sm_i \vu_i = \zeron\ \forall i \in [\![1,s]\!].
	\]
	Using again Lemma~\ref{lem:proppss}(ii), this latter property is equivalent 
	to $\set_i$ being a PSS for its linear span. Let $n_i$ be the 
	dimension of this span, and $m_i$ the number of elements in $\set_i$. 
	Since the $\Sm_i$ are orthogonal and the columns of $\Dm_{n,s}$ span 
	$\R^n$, we must have $\sum_{i=1}^s n_i=n$. In addition, we also have 
	$\sum_{i=1}^s m_i = n+s = \sum_{i=1}^s n_i + s$. Since $m_i >n_i$ for all 
	$i \in [\![1,s]\!]$, we conclude that $m_i=n_i+1$, and thus every $\set_i$ is 
	a minimal positive basis.
\end{proof}\\

Note that the characterization stated by Theorem~\ref{thm:ospbblocks} trivially 
holds for minimal positive bases. This theorem thus provides a 
characterization of the OSPB property through Gram matrices. One can then use 
this result to determine whether a positive basis has an orthogonal structure 
and, in that event, to highlight the associated decomposition.

\begin{corollary}
\label{cor:detOSPB}
	Let $n\ge 2$, $s \ge 1$ and $\D_{n,s}$ be a positive basis in $\R^n$. Let 
	$\Dm_{n,s}$ be a matrix representation of $\D_{n,s}$. If there exists a 
	permutation matrix $\Pm \in \R^{(n+s)\times (n+s)}$ such that the Gram 
	matrix $\G(\Dm_{n,s}\Pm)$ is block diagonal with $s$ blocks, then $\D_{n,s}$ 
	is an OSPB, and its decomposition~\eqref{eq:decompospb} is given by the 
	columns of $\Dm_{n,s}\Pm$ in that order.
\end{corollary}

Corollary~\ref{cor:detOSPB} provides a principled way of detecting the OSPB 
structure, by checking all possible permutations of the elements in the 
basis. Although this is not the focus of this work, we remark that it can be done in an efficient manner by considering a graph whose Laplacian matrix $\mathbf{L}$ has non-zero coordinates in the exact same rows and columns as $\G(\Dm_{n,s})$. Indeed, the problem of finding the 
permutation for matrix $\mathbf{L}$ reduces to that of finding the number of connected components in 
the graph~\cite{marsden2013eigenvalues} and can be solved in polynomial time~\cite{Hopcroft1973algograph}.

%%%%%%%%%%%%%%%%%%%%%%%%%%%%%%%%%%%%%%%%%%%%%%%%%%%%%%%%%%%%%%%%%%%%%%%%%%%%%%%
\subsection{Cosine measure of an orthogonally structured positive basis}
\label{subsec:cmOSPB}

As seen in the previous section, orthogonally structured positive bases admit 
a decomposition into minimal positive bases that are orthogonal to one 
another. In this section, we investigate how this particular structure allows 
for efficient computation of the cosine measure.

Our approach builds on the algorithm introduced by Hare and 
Jarry-Bolduc~\cite{hare2020cm}, that computes the cosine measure of any 
positive spanning set in finite time (though the computation time may be 
exponential in the problem dimension). This procedure consists of a loop 
over all possible linear bases contained in the positive spanning set. More 
formally, given a positive spanning set $\D$ in $\R^n$ and a linear basis 
$\B \subset \D$, we let $\Dm$ be a matrix representation of $\D$, and $\Bm$ 
the associated representation of $\B$. One then defines
\begin{equation}
\label{eq:qtiesBcmalgo}
	\gamma_{\B} := \frac{1}{\sqrt{\onen^\T \G(\Bm)^{-1}\onen}} 
	\quad \mbox{and} \quad
	\vu_{\B} := \gamma_{\B} \Bm^{-\T}\onen.
\end{equation}
The vector $\vu_{\B}$ is the only unitary vector such that 
$\vu_{\B}^\T \tfrac{\vd}{\|\vd\|} = \gamma_{\B}$ for all 
$\vd \in \B$~\cite[Lemmas 12 and 13]{hare2020cm}. Computing the 
quantities~\eqref{eq:qtiesBcmalgo} for all linear bases contained in $\D$ 
then gives both the cosine measure and the cosine vector set for 
$\D$~\cite[Theorem 19]{hare2020cm}. Indeed, we have
\begin{equation}
\label{eq:cmformulaHJB}
	\cm{\D}=\min_{\substack{\B \, \subset \, \D \\ \B\ basis\ of\ \R^n}} 
	\max_{\vd \in \D} \vu_{\B}^\T \frac{\vd}{\|\vd\|},
\end{equation}
while the cosine vector set is given by
\[
	\cV(\D)= \left\{ 
	\vu_{\B}:\max_{\vd \in \D} \vu_{\B}^\T \tfrac{\vd}{\|\vd\|} = \cm{\D}
	\right\}.
\]
The next proposition shows how, when dealing with OSPBs, 
formula~\eqref{eq:cmformulaHJB} can be simplified to give a more direct link 
between the quantities $\gamma_{\B}$ and the cosine measure.

\begin{proposition}
\label{prop:OSPBcmgamma}
	Let $\D_{n,s}$ be an OSPB of $\R^n$ and $\D_{\LL_1},\dots,\D_{\LL_s}$ be a 
	decomposition of $\D_{n,s}$ into minimal positive bases. Then,
	$$
		\cm{\D_{n,s}}=\min_{\substack{\B_n \, \subset \, \D_{n,s} \\ 
		\B_n\ basis\ of\ \R^n}}
		\gamma_{\B_n} 
		\quad \mbox{and} \quad 
		\cV(\D_{n,s})=\{\vu_{\B_n}:\gamma_{\B_n}=\cm{\D_{n,s}}\},
	$$
	where $\gamma_{\B_n}$ and $\vu_{\B_n}$ are computed according 
	to~\eqref{eq:qtiesBcmalgo}. 
\end{proposition}

\begin{proof}
	Using the decomposition of $\D_{n,s}$, any linear basis $\B_n$ of $\R^n$ contained in $\D_{n,s}$ can 
	be decomposed as $\B_n=\B_{\LL_1} \cup \cdots \cup \B_{\LL_s}$ where 
	$\B_{\LL_i} \subset \D_{\LL_i}$ is a linear basis of $\LL_i$. By 
	construction, the vector $\vu_{\B_n}$ makes a positive dot product 
	with every element of $\B_n$ (and thus with every element of any 
	$\B_{\LL_i}$) equal to $\gamma_{\B_n}$. Consequently, the vector $\vu_{\B_n}$ 
	lies in the positive span of $\B_n$, and its projection on any subspace 
	$\LL_i$ lies in the positive span of $\B_{\LL_i}$.
	
 	Meanwhile, for any $i \in [\![1,s]\!]$, there exists $\vd_i \in \D_{\LL_i}$ 
	such that $\B_{\LL_i}=\D_{\LL_i}\setminus \{\vd_i\}$, and this vector 
	satisfies	$\vu_{\B_n}^\T \vd_i < 0$ per 
	Lemma~\ref{lem:proppss}$(\ref{lem: PSS eq pos dotprod})$, leading to
	\begin{equation*}
		\max_{\vd \in \D_n} \frac{\vu_{\B_n}^\T \vd}{\|\vd\|} 
		= \max_{i \in [\![1,s]\!]} \left\{ 
		\max_{\vd \in \D_{\LL_i}} \frac{\vu_{\B_n}^\T \vd}{\|\vd\|} \right\} 
		= \max_{i \in [\![1,s]\!]} \left\{
		\max_{\vd \in \B_{\LL_i}} \frac{\vu_{\B_n}^\T \vd}{\|\vd\|} \right\} 
		= \max_{\vd \in \B_n} \frac{\vu_{\B_n}^\T \vd}{\|\vd\|} =\gamma_{\B_n},
	\end{equation*}
	where the second equality comes from $\vu_{\B_n}^\T \vd_i < 0$ and 
	$\max\limits_{\vd \in \D_n} \vu_{\B_n}^\T \vd >0$, and the last equality 
	holds by definition of $\gamma_{\B_n}$. Recalling that $\cm{\D_n}$ is given 
	by~\eqref{eq:cmformulaHJB} concludes the proof.
\end{proof}\\

One drawback of the approach described so far is that it requires 
computing the quantities~\eqref{eq:qtiesBcmalgo} for all linear bases including 
in the positive spanning set of interest. In the case of OSPBs, we show below 
that the number of linear bases to be considered can be greatly reduced by 
leveraging to their decomposition into minimal positive bases.

\begin{theorem} 
\label{thm:OSPBcm}
	Let $\D_{n,s}$ be an OSPB of $\R^n$ and $\D_{\LL_1},\dots,\D_{\LL_s}$ be a 
	decomposition of $\D_{n,s}$ into minimal positive bases. Let $\B_n$ be a 
	linear basis contained in $\D_{n,s}$ such that 
	$\B_n = \cup_{i=1}^s \B_{\LL_i}$, where every $\B_{\LL_i}$ is a linear 
	basis of the subspace corresponding to $\D_{\LL_i}$. For each basis 
	$\B_{\LL_i}$, define
	$$
		\gamma_{\B_{\LL_i}} = \frac{1}{
		\sqrt{\one_{\LL_i}^\T \G(\Bm_{\LL_i})^{-1} \one_{\LL_i}}},
	$$
	where $\Bm_{\LL_i}$ is a matrix representation of $\B_{\LL_i}$ and 
	$\one_{\LL_i}$ denotes the vector of all ones in $\R^{\dim(\LL_i)}$. Then,
	\begin{equation}
	\label{eq:OSPBgamma}
		\gamma_{\B_n} = \frac{1}{\sqrt{\sum_{i=1}^s \gamma_{\B_{\LL_i}}^{-2}}}
		= \frac{1}{\sqrt{\sum_{i=1}^s \one_{\LL_i}^\T \G(\Bm_{\LL_i})^{-1} 
		\one_{\LL_i}}},
	\end{equation}
	with $\gamma_{\B_n}$ being defined as in~\eqref{eq:qtiesBcmalgo}.
	
	As a result, the cosine measure and cosine vector set of $\D_{n,s}$ are 
	given by
	\begin{equation}
	\label{eq:OSPBcm}
		\cm{\D_{n,s}} = \frac{1}{\sqrt{\sum_{i=1}^s
		\max_{\vd_i \in \D_{\LL_i}} 
		\gamma_{\D_{\LL_i} \setminus \{\vd_i\}}^{-2}}}
	\end{equation}
	and
	\begin{equation}
	\label{eq:OSPBcv}
		\cV(\D_{n,s}) = \left\{\vu_{\B_n}: \B_n = \cup_{i=1}^s \B_{\LL_i}, 
		\quad \vu_{\B_n} = \frac{ 
		[\Bm_{\LL_1}\ \cdots\ \Bm_{\LL_s}]^{-\T}\one_n}
		{\sqrt{\sum_{i=1}^s \gamma_{\B_{\LL_i}}^{-2}}}
		\right\}
	\end{equation}
	respectively.
\end{theorem}

\begin{proof}
	Let $\Bm_n$ be a matrix representation of $\B_n$ such that 
	$\Bm_n = \bbm \Bm_{\LL_1} \cdots \Bm_{\LL_s} \ebm$. Then, the Gram matrix 
	of $\Bm_n$ is 
	$\G(\Bm_n)=\diag\left( \G(\Bm_{\LL_1}),\dots,\G(\Bm_{\LL_s})\right)$, 
	implying that
	\begin{eqnarray*}
		\gamma_{\B_n} 
		= \frac{1}{\sqrt{\onen^\T \G(\Bm_n)^{-1} \onen}} 
		&= &\frac{1}{\sqrt{\onen^\T \diag\left( 
		\G(\Bm_{\LL_1}), \dots, \G(\Bm_{\LL_s})\right )^{-1} \onen}} \\
		&= &\frac{1}{\sqrt{\onen^\T \diag\left( 
		\G(\Bm_{\LL_1})^{-1}, \dots, \G(\Bm_{\LL_s})^{-1}\right) \onen}} \\
		&= &\frac{1}{\sqrt{\sum_{i=1}^s 
		\one_{\LL_i}^\T \G(\Bm_{\LL_i})^{-1}\one_{\LL_i}}},
	\end{eqnarray*}
	which proves~\eqref{eq:OSPBgamma}. 
	
	Recall now that every linear basis $\B_{\LL_i}$ can be written as 
	$\D_{\LL_i} \setminus \{\vd_i\}$ for some $\vd_i \in \D_{\LL_i}$.
	Combining this property together with~\eqref{eq:OSPBgamma} and 
	Proposition~\ref{prop:OSPBcmgamma}, we obtain
	\begin{eqnarray*}
		\cm{\D_{n,s}} 
		= \min_{\substack{\B_n \subset \D_{n,s} \\
		\B_n\ basis\ of\ \R^n}}\gamma_{\B_n} 
		&= &\min_{\substack{\B_{\LL_1},\dots,\B_{\LL_s} \\
		\B_{\LL_i} = \D_{\LL_i} \setminus \{\vd_i\}}} 
		\frac{1}{\sqrt{\sum_{i=1}^s \gamma_{\B_{\LL_i}}^{-2}}} \\
		&= &\min_{\substack{\vd_1,\dots,\vd_s \\
		\vd_i \in \D_{\LL_i}}} 
		\frac{1}{\sqrt{\sum_{i=1}^s 
		\gamma_{\D_{\LL_i} \setminus \{\vd_i\}}^{-2}}} \\
		&= &\frac{1}{\sqrt{\sum_{i=1}^s
		\max_{\vd_i \in \D_{\LL_i}} 
		\gamma_{\D_{\LL_i} \setminus \{\vd_i\}}^{-2}}},
	\end{eqnarray*}
	proving~\eqref{eq:OSPBcm}. 
	
	Finally, combining~\eqref{eq:OSPBcm} with the result of 
	Proposition~\ref{prop:OSPBcmgamma} on the cosine vector set gives
	\begin{eqnarray*}
		\cV(\D_{n,s})
		&= &\left\{\vu_{\B_n}:\gamma_{\B_n}=\cm{\D_{n,s}}\right\} \\
		&= &\left\{\vu_{\B_n}: 
		\B_n = \cup_{i=1}^s \B_{\LL_i},\quad \frac{1}{\sqrt{\sum_{i=1}^s
		\gamma_{\B_{\LL_i}}^{-2}}} = \cm{\D_{n,s}} \right\} 
	\end{eqnarray*}
	Using a matrix representation $\Bm_n=\bbm \Bm_{\LL_1} \cdots \Bm_{\LL_s} 
	\ebm$ along with the formula~\eqref{eq:qtiesBcmalgo} for $\vu_{\B_n}$ 
	leads to
	\begin{eqnarray*}
		\cV(\D_{n,s})
		&= &\left\{\vu_{\B_n}: \B_n = \cup_{i=1}^s \B_{\LL_i},
		\quad \vu_{\B_n} = \gamma_{\B_n} 
		[\Bm_{\LL_1}\ \cdots\ \Bm_{\LL_s}]^{-\T}\one_n
		\right\} \\
		&= &\left\{\vu_{\B_n}: \B_n = \cup_{i=1}^s \B_{\LL_i}, 
		\quad \vu_{\B_n} = \frac{ 
		[\Bm_{\LL_1}\ \cdots\ \Bm_{\LL_s}]^{-\T}\one_n}
		{\sqrt{\sum_{i=1}^s \gamma_{\B_{\LL_i}}^{-2}}}
		\right\} \\
	\end{eqnarray*}
	which is precisely~\eqref{eq:OSPBcv}.
\end{proof}\\

On a broader level, Theorem~\ref{thm:OSPBcm} shows that any calculation for 
a linear basis contained in an OSPB reduces to calculations on bases contained 
in each of the minimal positive bases in its decomposition. This observation 
leads to a principled way of computing the cosine measure from the OSPB 
decomposition, as described in Algorithm~\ref{algo:cmOSPB}.

%%%%%%%%%%%%%%%%%%%%%%%%%%%%%%%%%%%%%%%
\noindent 
\begin{algorithm} [H] \label{algo:cmOSPB}
	\DontPrintSemicolon
	\caption{Cosine measure of an orthogonally structured positive basis}

	\emph{Input:} An orthogonally structured positive basis $\D_{n,s}$ of 
	$\R^n$ with normalized vectors.

	\vspace{0.2cm}

	\textbf{Step 1.} Compute a decomposition 
	$\D_{n,s} = \D_{\LL_1} \cup \cdots \cup \D_{\LL_s}$ of $\D_{n,s}$ into $s$ 
	minimal positive bases $\D_{\LL_i}$ on subspaces of $\R^n$. 

	\vspace{0.2cm}

	%\textbf{Step 2.} For every $i \in [\![1,s]\!]$, compute
	%$$
	%	\beta_{\LL_i} = \max_{\vd \in \D_{\LL_i}} 
	%	\gamma_{\D_{\LL_i} \setminus \{\vd\}}^{-2}
	%	\quad \mbox{and} \quad
	%	\B_{\LL_i} \in \argmax{\vd \in %\D_{\LL_i}}\ 
	%	\gamma_{\D_{\LL_i} \setminus \{\vd\}}^{-2}.
	%$$
\textbf{Step 2.} For every $i \in [\![1,s]\!]$, compute
	$$
		\beta_{\LL_i} = \max_{\vd \in \D_{\LL_i}} 
		\gamma_{\D_{\LL_i} \setminus \{\vd\}}^{-2}
		\quad \mbox{and} \quad
		\mathcal{A}_i=\argmax{\vd \in \D_{\LL_i}}\ 
		\gamma_{\D_{\LL_i} \setminus \{\vd\}}^{-2}.
	$$ For all $\vd \in \mathcal{A}_i,$ note $\B_{\LL_i}^{(\vd)}=\D_{\LL_i}\setminus\{\vd\}.$ 

	\vspace{0.2cm}

	\textbf{Step 3.} Return the cosine measure
	\begin{equation}
	\label{eq:cmOSPBalgo}
		\cm{\D_{n,s}} = \frac{1}{\sqrt{\sum_{i=1}^s \beta_{\LL_i}}}
	\end{equation}
	and the cosine vector set
	$$
		\cV(\D_{n,s}) = \bigcup\limits_{(\vd_1,\dots,\vd_s)\in \mathcal{A}_1\times\cdots\times\mathcal{A}_s}\left\{
		 \vu \in \R^n, \vu = 
		\frac{\bbm \Bm_{\LL_1}^{(\vd_1)}&\cdots&\Bm_{\LL_s}^{(\vd_s)} \ebm ^{-\T} \onen}
		{\sqrt{\sum_{i=1}^s \beta_{\LL_i}}} 
		\right\},
	$$
	where $\Bm_{\LL_i}^{(\vd)}$ is a matrix representation of $\B_{\LL_i}^{(\vd)}.$
\end{algorithm}
%%%%%%%%%%%%%%%%%%%%%%%%%%%%%%%%%%%%%%%

\vspace{0.2cm}
%\noindent 
In essence, Algorithm~\ref{algo:cmOSPB} is quite similar to the generic 
method proposed for positive spanning sets~\cite{hare2020cm}. However, 
Algorithm~\ref{algo:cmOSPB} reduces the computation of the cosine measure 
to that over minimal positive bases, which leads to significantly cheaper 
calculations. Indeed, for any minimal positive basis $\D_{\LL_i}$ involved 
in the decomposition, the algorithm considers $|\D_{\LL_i}|$ positive bases, 
and computes the quantities of interest~\eqref{eq:qtiesBcmalgo} for each 
of those (which amounts to inverting or solving a linear system involving the 
associated Gram matrix). Overall, Algorithm~\ref{algo:cmOSPB} thus 
checks $\sum_{i=1}^s |\D_{\LL_i}|=|\D_{n,s}|=n+s $ bases to 
compute the cosine measure~\eqref{eq:cmOSPBalgo}. This result 
represents a significant improvement over the $\binom{n+s}{n}$ possible 
bases that are potentially required when the OSPB decomposition is not 
exploited, as in the algorithm for generic PSSs~\cite{hare2020cm}. We also 
recall from Section~\ref{subsec:structOSPB} that Step 1 in 
Algorithm~\ref{algo:cmOSPB} can be performed in polynomial time.

To end this section, we illustrate how Algorithm~\ref{algo:cmOSPB} results in 
a straightforward calculation in the case of some minimal positive bases, by 
computing explicitly the cosine measure of the minimal positive basis 
$\In \cup \{-\onen\}$. In this case, the calculation becomes easy thanks to the orthogonality of the vectors within $\In \cup \{-\one\}$. Although the value~\eqref{eq:cmIn1n} was recently 
stated~\cite{jarrybolducthesis,roberts2022subspace} and checked numerically 
using the method of Hare and Jarry-Bolduc~\cite{hare2020cm}, to the best of 
our knowledge the formal proof below is new. 

\begin{lemma}
\label{lem:cmIn1n}
	The cosine measure of the OSPB $\In \cup \{-\onen\}$ is given by
	\begin{equation}
	\label{eq:cmIn1n}
		\cm{\In \cup \{-\onen\}} = \frac{1}{\sqrt{n^2+2(n-1)\sqrt{n}}}.
	\end{equation}
\end{lemma}
\begin{proof}
	For the sake of simplicity, we normalize the last vector in the set (which does 
	not change the value of the cosine measure) and we let 
	$\D_n=\mathcal{I}_n \cup\{-\tfrac{1}{\sqrt{n}}\onen \}$ in the 
	rest of the proof. Since $\D_n$ is a positive basis, it follows that
	\[
		\cm{\D_n}=\min_{\vd \in \D_n} \gamma_{\D_n \setminus \{\vd\}}, 
	\] 
	Two cases are to be considered. Suppose first that 
	$\vd=-\tfrac{1}{\sqrt{n}}\onen$, so that $\D_n \setminus \{\vd\}=\In$. Then, 
	any matrix representation of that basis $\Bm$ is such that $\G(\Bm)=\mIn$. 
	Consequently,
	\[
		\gamma_{\D_n \setminus \{\vd\}} 
		= \frac{1}{\sqrt{\onen^\T \G(\Bm)^{-1}\onen}} 
		= \frac{1}{\sqrt{\onen^\T \onen}} = \frac{1}{\sqrt{n}}.
	\]
	Suppose now that $\vd=\ve_i$ for some $i \in [\![1,n]\!]$. In 
	that case, considering the matrix representation 
	$\Bm_i = \bbm \ve_1 &\cdots &\ve_{i-1} &\ve_{i+1} &\cdots &\ve_n 
	&-\tfrac{1}{\sqrt{n}}\onen \ebm$, we obtain that 
	\[
		\G(\Bm_i) = \bbm \mI_{n-1} &-\tfrac{1}{\sqrt{n}}\one_{n-1} \\ 
		-\tfrac{1}{\sqrt{n}}\one_{n-1}^\T &1 \ebm
		\quad \mbox{and} \quad
		\G(\Bm_i)^{-1} = \bbm \mI_{n-1}+\one_{n-1}\one_{n-1}^\T &\sqrt{n}\one_{n-1}\\ 
		\sqrt{n}\one_{n-1}^\T &n \ebm
	\]
	Since these formulas do not depend on $i$, we obtain that for all 
	$i \in [\![1,n]\!]$, we have
	\[
		\gamma_{\D_n \setminus \{\ve_i\}} 
		= \frac{1}{\sqrt{\onen^\T \G(\Bm_i)^{-1} \onen}} 
		= \frac{1}{\sqrt{\sum_{j=1}^n \sum_{\ell=1}^n [\G(\Bm_i)^{-1}]_{j\ell}}}.
	\]
	Summing all coefficients of $\G(\Bm_i)^{-1}$ gives 
	$n-1+(n-1)^2 + 2(n-1)\sqrt{n}+n=n^2 +2(n-1)\sqrt{n}$, hence 
	$\gamma_{\D_n \setminus \{\ve_i\}} = \frac{1}{\sqrt{n^2+2(n-1)\sqrt{n}}}$.
	Comparing this value with $\tfrac{1}{\sqrt{n}}$ yields the desired conclusion.
%	\crnote{The proof does not seem to work when the last vector is not normalized.}
\end{proof}\\

Algorithm~\ref{algo:cmOSPB} allows for efficient calculation of cosine 
measures for specific positive spanning sets that originate from OSPBs, as we 
will establish in Section~\ref{sec:pksscmk} in the case of positive 
$k$-spanning sets. 

%%%%%%%%%%%%%%%%%%%%%%%%%%%%%%%%%%%%%%%%%%%%%%%%%%%%%%%%%%%%%%%%%%%%%%%%%%%%%%%
\section{Positive \texorpdfstring{$k$}{TEXT}-spanning sets and \texorpdfstring{$k$}{TEXT}-cosine measure}
\label{sec:pksscmk}
%%%%%%%%%%%%%%%%%%%%%%%%%%%%%%%%%%%%%%%%%%%%%%%%%%%%%%%%%%%%%%%%%%%%%%%%%%%%%%%

In this section, we are interested in positive spanning sets that retain their 
positively spanning ability when one or more of their elements are removed. 
Those sets were originally termed \emph{positive 
$k$-spanning sets}~\cite{Marcus1984}, and we adopt the same terminology. 
Section~\ref{subsec:pkssdefs} recalls the key definitions for positive 
$k$-spanning sets and positive $k$-bases, while Section~\ref{subsec:cmk} 
introduces the $k$-cosine measure, a generalization of the cosine measure from 
Section~\ref{subsec:cm}. Finally, Section~\ref{subsec:pkssospb} illustrates 
how to construct sets with guaranteed $k$-cosine measure based on OSPBs.

%%%%%%%%%%%%%%%%%%%%%%%%%%%%%%%%%%%%%%%%%%%%%%%%%%%%%%%%%%%%%%%%%%%%%%%%%%%%%%%
\subsection{Positive \texorpdfstring{$k$}{TEXT}-spanning property}
\label{subsec:pkssdefs}

Our goal for this section is to provide a summary of results on positive 
$k$-spanning sets that mimic the standard ones for positive spanning sets. We 
start by defining the property of interest.

\begin{definition}[Positive $k$-span and positive $k$-spanning set]
\label{def:pkspanpkss}
	Let $\LL$ be a linear subspace of $\R^n$ and $k \ge 1$. The 
	\emph{positive $k$-span} of a finite family of vectors 
	$\D$ in $\LL$, 
	denoted $\pspan_k(\D)$, is the set 
	$$
		\pspan_k(\D):=
		\bigcap_{\substack{\setS \subset \D \\ 
		|\setS|\geq|\mathcal{D}|-k+1}}\pspan(\setS).
	$$
	A \emph{positive $k$-spanning set} (P$k$SS) of $\LL$ is a family of vectors 
	$\D$ such that $\pspan_k(\D)=\LL$.
\end{definition}

As for Definition~\ref{def:pspanpss}, when $\LL$ is not specified, a P$k$SS 
should be understood as a P$k$SS of $\R^n$. By construction, any P$k$SS of 
$\LL$ must contain a PSS of $\LL$, and therefore is a PSS of $\LL$ itself. 
Moreover, the notions of positive $k$-span and P$k$SS with $k=1$ coincide with 
that of positive span and PSS from Definition~\ref{def:pspanpss}. Similar to 
PSSs, we will omit the subspace of interest when $\LL=\R^n$.

The notion of positive spanning set is inherently connected to that of 
spanning set, a standard concept in linear algebra. We provide below a 
counterpart notion associated with P$k$SSs.

\begin{definition}[$k$-span and $k$-spanning set]
\label{def:kspan}
	Let $\D$ be a finite family of vectors in $\R^n$. the \emph{$k$-span} of 
	$\D$, denoted $\spann_k(\D)$, is defined by
	\[
		\spann_k(\D) \; = \; \bigcap_{\substack{\setS \subset \D \\ 
		|\setS|\geq|\D|-k+1}} \spann(\setS).
	\]
	Given a subspace $\LL$ of $\R^n$, a \emph{$k$-spanning set} of $\LL$ is a 
	family of vectors $\D$ such that $\spann_k(\D)=\LL$.
\end{definition}

Using the two definitions above, we can generalize the results of 
Lemma~\ref{lem:proppss} and Lemma \ref{lem:pss => linspanning} to positive 
$k$-spanning sets. The proof is omitted as it merely amounts to combining the 
definitions with these two lemmas.

\begin{lemma}
\label{lem:proppkss}
	Let $\LL$ be a subspace of $\R^n$ and $\D$ a finite family of vectors in 
	$\LL$. Let $k \in \mathbb{N}$ satisfy $1\leq k \leq |\D|$. The following statements are equivalent:
	\begin{enumerate}[(i)]
        \item\label{equiv D pkss} $\D$ is a P$k$SS of $\LL$.
	    \item\label{equiv k pos dot prod} For any nonzero vector 
	    $\vu \in \LL$, there exist $k$ elements $\vd_1,\dots,\vd_k$ of $\D$ such 
	    that $\vu^\T \vd_i>0$ for all $i \in [\![1,k]\!]$.
    	\item\label{equiv kpsan and 0pos} $\spann_k(\D)=\LL$ and 
  		for any $\setS \subset \D$ of cardinality $|\D|-k+1$, the vector $\zeron$ 
    	can be written as a \textbf{positive} linear combination of the elements 
    	of $\setS$.
 	\end{enumerate}
\end{lemma}

\begin{lemma}
\label{lem:pkss => linspanning}
   Let $\LL$ be a subspace of $\R^n$ and let the finite set $\D$ be a P$k$SS for 
   $\LL$. Then, for any $\vd \in \D$, the $k$-span of the family 
   $\D \setminus \{\vd\}$ is $\LL$.
\end{lemma}

The equivalence between statements $(\ref{equiv D pkss})$ and 
$(\ref{equiv k pos dot prod})$ from Lemma~\ref{lem:proppkss} motivates the 
term ``positive $k$-spanning sets''. Indeed, given a P$k$SS and a vector in the 
subspace it positively $k$-spans, there exist $k$ elements of the P$k$SS that 
make an acute angle with that vector. Note however that the equivalence between 
statements $(\ref{equiv D pkss})$ and $(\ref{equiv kpsan and 0pos})$, as well 
as Definition~\ref{def:pkspanpkss}, both provide an alternate characterization, 
namely that a P$k$SS is a PSS that retains this property when removing $k-1$ 
of its elements. This latter characterization implies that a P$k$SS must 
contain at least $n+k$ vectors, although this bound is only tight when $k=1$. 
Indeed, in \cite[Corollary 5]{Marcus1984} it is shown that the minimal size of 
a positive $k$-spanning set in a subspace of dimension $\ell$ is 
$2k+\ell-1$.

Throughout Section~\ref{sec:defs}, we highlighted positive bases as a special 
class of positive spanning sets. We now define the counterpart notion for 
positive $k$-spanning sets, that are termed positive $k$-bases.

\begin{definition}[Positive $k$-basis]\label{def:pkb}
	Let $\LL$ be a linear subspace of $\R^n$. A \emph{positive $k$-basis} of 
	$\LL$ is a P$k$SS $\D$ of $\LL$ satisfying
  	\[
  		\forall \mathbf{d} \in \D, \quad
  		\pspan_k{(\D\backslash\{\mathbf{d}\})} \neq \LL.
	\]
\end{definition}

Positive $k$-bases can be thought as inclusion-wise minimal positive 
$k$-spanning sets (similarly to the characterization of positive bases from 
Section~\ref{subsec:pss}). We showed earlier how the notion of positive 
independence can be used to give an alternative definition for positive bases. 
Let us generalize this idea and introduce the concept of positive 
$k$-independence, used to characterize positive $k$-bases.
\begin{definition}[Positive $k$-independence]
\label{def:pki}
 	Let $\LL$ be a linear subspace of $\R^n$ of dimension $\ell \ge 1$ and 
	let $k \ge 1$. A family of vectors $\D$ in $\LL$ is \emph{positively 
	$k$-independent} if $|\D| \geq k$ and, for any $\vd \in \D$, there exists a vector 
	$\vu \in \LL$ having a positive dot product with \textbf{exactly} $k$ 
	elements of $\D$, including $\vd$.
\end{definition}
Using this new concept, positive $k$-bases can alternatively be defined as 
positively $k$-independent P$k$SSs. We mention in passing that the upper bound 
on the size of a positive $k$-basis has yet to be determined, though it is 
known to exceed $2k\ell$ for an $\ell$-dimensional 
subspace~\cite{Wotzlaw2009}.

%%%%%%%%%%%%%%%%%%%%%%%%%%%%%%%%%%%%%%%%%%%%%%%%%%%%%%%%%%%%%%%%%%%%%%%%%%%%%%%
\subsection{The \texorpdfstring{$k$}{TEXT}-cosine measure}
\label{subsec:cmk}

This section aims at generalizing the cosine measure from 
Section~\ref{subsec:cm} to positive $k$-spanning sets so that it characterizes 
the positive $k$-spanning property. Our proposal, called the $k$-cosine 
measure, is described in Definition~\ref{def:cmk}.

\begin{definition}[$k$-cosine measure]
\label{def:cmk}
	Let $\D$ be a finite family of nonzero vectors in $\R^n$ and let $k \in \mathbb{N}$ satisfy $1\leq k \leq |\D|$. 
	The \emph{$k$-cosine measure} of $\D$ is given by
	$$
		\cmk{\D}{k} = \min_{\substack{\|\vu\| =1 \\ \vu \in \R^n}}
		\max_{\substack{\setS \subset \D \\ |\setS|=k}} 
		\min_{\vd \in \setS} \frac{\vu^\T \vd}{\|\vd\|}.
	$$
	The \emph{$k$-cosine vector set} associated with $\D$ is given by
	\[
		\cV_k(\D) = \argmin{\substack{\|\vu\| =1 \\ \vu \in \R^n}} 
		\max_{\substack{\setS \subset \D \\ |\setS|=k}} 
		\min_{\vd \in \setS} \frac{\vu^\T \vd}{\|\vd\|}.
	\]
\end{definition}

Note that Definition~\ref{def:cm} is a special case of 
Definition~\ref{def:cmk} corresponding to $k=1$. In its general form, this 
definition expresses how well the vectors of the family of interest are 
spread in the space through subsets of $k$ vectors, which is related to 
property $(\ref{equiv k pos dot prod})$ in Lemma~\ref{lem:proppkss}. Our next result shows that 
this quantity characterizes P$k$SSs.

\begin{theorem} 
\label{thm:cmkpkss}
	Let $\D$ be a finite family of vectors in $\R^n$ and let $k \in \mathbb{N}$ satisfy $1\leq k \leq |\D|$. Then $\D$ is a 
	positive $k$-spanning set in $\R^n$ if and only if $\cmk{\D}{k}>0$.
\end{theorem}
\begin{proof}
	Without loss of generality we assume that all the elements of $\D$ are unit 
	vectors. 
	
	Suppose first that $\D$ is a $P$k$SS$. Then, by Lemma~\ref{lem:proppkss}, 
	for any unit vector $\vu \in \R^n$, there exist $k$ vectors in $\D$ that 
	make a positive dot product with $\vu$. Let $\set_{\vu}$ be a subset of 
	$\D$ consisting of these $k$ vectors. By construction, we have 
	$\min\limits_{\vd \in \set_{\vu}} \vu^\T \vd>0$ and thus
	\[
		\max_{\substack{\setS \subset \D \\ |\setS|=k}} 
		\min_{\vd \in \setS} \vu^\T \vd
		\ge 
		\min_{\vd \in \set_{\vu}} \vu^\T \vd >0.
	\]
	Since the result holds for any unit vector $\vu$, we conclude that 
	$\cmk{\D}{k}>0$.
	
	Conversely, suppose that $\cmk{\D}{k}>0$. For any unit vector $\vu$, we have by assumption 
	that
	\begin{equation}
	\label{eq:cmkpksscvsly}
		\max_{\substack{\setS \subset \D \\ |\setS|=k}} 
		\min_{\vd \in \setS} \vu^\T \vd >0.
	\end{equation} 
In other words, at least $k$ elements of $\D$ have a positive dot product with $\vu$. This proves that $\D$ satisfies statement~$(\ref{equiv k pos dot prod})$ from 
	Lemma~\ref{lem:proppkss} and thus $\D$ is a positive $k$-spanning set.
 
	%Let $\bar{\setS}_1$ be a subset of $\D$ such that $|\bar{\set}_1|=|\D|-k+1$.
 	%Then, \eqref{eq:cmkpksscvsly} implies that
	%\[
		%\max_{\substack{\setS \subset \bar{\set}_1 \\ |\setS|=1}} 
		%\min_{\vd \in \setS} \vu^\T \vd 
		%= 
		%\max_{\vd \in \bar{\set}_1} \vu^\T \vd >0,
	%\]
	%hence there exists a vector $\bar{\vd}_1 \in \bar{\set}_1$ such that 
	%$\vu^\T \bar{\vd}_1>0$. Consider now a set $\bar{\set}_2$ such that 
	%$\bar{\set}_2=|\D|-k+1$ and $\bar{\vd}_1 \notin \bar{\set}_2$. Applying the 
	%same reasoning as above shows that there must exist 
	%$\bar{\vd}_2 \in \bar{\set}_2$ such that $\vu^\T \bar{\vd}_2>0$. Overall, 
	%we can repeat the process $k$ times, showing that there must exist at least 
	%$k$ elements of $\D$ that make a positive dot product with $\vu$. Since the 
	%result applies to any unit vector, we have shown 
	%Lemma~\ref{lem:proppkss}$(\ref{equiv k pos dot prod})$, and thus that $\D$ is 
	%a P$k$SS.
\end{proof}\\

As announced in Section~\ref{subsec:cm}, the proof of 
Theorem~\ref{thm:cmkpkss} covers that of Proposition~\ref{prop:cmpss} as the 
special case $k=1$.

To end this section, we provide an alternate definition of the 
$k$-cosine measure that connects this notion to the cosine measure.

\begin{theorem}
\label{thm:othercmk}
	Let $\D$ be a finite family of nonzero vectors in $\R^n$ and let $k \in \mathbb{N}$ satisfy $1\leq k \leq |\D|$. Then,
	\begin{equation}
	\label{eq:othercmk}
		\cmk{\D}{k} \; = \; 
		\min_{\substack{\setS \subset \D  \\ |\setS|=|\D|-k+1}} \cm{\setS}.
	\end{equation}
\end{theorem}
\begin{proof}
	Without loss of generality, we assume that 
	$\D=\{\vd_1,\dots,\vd_m\}$ where each $\vd_i$ is a unit vector. 
	In order to show that
	\[
		\min\limits_{\|\vu\|=1}
		\max_{\substack{\set \subset \D \\ |\set|=k}}
		\min\limits_{\vd \in \set} \vu^{\T}\vd
		=
		\min_{\substack{\set \subset \D \\|\set|=m-k+1}}
		\cm{\set},
	\]
	we prove the equivalent statement
	\begin{equation}
	\label{eq:othercmkequiv}
		\min\limits_{\|\vu\|=1}
		\max_{\substack{\set \subset \D \\ |\set|=k}}
		\min\limits_{\vd \in \set} \vu^{\T}\vd 
		=
		\min\limits_{\|\vu\|=1}
		\min_{\substack{\set \subset \D \\ |\set|=m-k+1}}
		\max\limits_{\vd \in \set} \vu^{\T}\vd.
	\end{equation}
	To this end, let $\vu$ be a unit vector in $\R^n$. We reorder the 
	elements of $\D$ so that for all $i \leq m-1$, $\vu^\T \vd_i \ge \vu^\T \vd_{i+1}$ and we define 
	$\set_k^{-}=\{\vd_1,\dots,\vd_k\}$ and $\set_k^+ = \{\vd_k,\dots,\vd_m\}$. 
	By construction, one has 
	\begin{equation}
	\label{eq:othercmkdktwosets}
		\vu^\T \vd_k 
		= \min_{\vd \in \set_k^{-}} \vu^\T \vd 
		= \max_{\vd \in \set_k^{+}} \vu^\T \vd.
	\end{equation}
	Moreover, the definitions of $\set_k^{-}$ and $\set_k^+$ also imply that
	\begin{equation}
	\label{eq:othercmdkargminargmax}
		\set_k^{-} \in 
		\argmax{\substack{\set \subset \D \\ |\set|=k}}
		\min_{\vd \in \set}\vu^{\T}\vd 
		\quad \mbox{and} \quad		
		\set_k^+ \in 
		\argmin{\substack{\set \subset \D \\ |\set|=m-k+1}} 
		\max_{\vd \in \set} \vu^\T \vd.
	\end{equation}
	Combining~\eqref{eq:othercmkdktwosets} and~\eqref{eq:othercmdkargminargmax} 
	gives
   	\begin{equation}
   	\label{eq:othercmkequivu}
   		\max_{\substack{\set \subset \D \\ |\set|=k}}
		\min_{\vd \in \set}\vu^{\T}\vd 
		\; = \;
		\min_{\substack{\set \subset \D \\ |\set|=m-k+1}} 
		\max_{\vd \in \set} \vu^\T \vd.
   	\end{equation}
	Since~\eqref{eq:othercmkequivu} holds for any unit vector $\vu \in \R^n$, it 
	follows that
	\[
		\min_{\|\vu\|=1}
   		\max_{\substack{\set \subset \D \\ |\set|=k}}
		\min_{\vd \in \set}\vu^{\T}\vd 
		\; = \;
		\min_{\|\vu\|=1}
		\min_{\substack{\set \subset \D \\ |\set|=m-k+1}} 
		\max_{\vd \in \set} \vu^\T \vd,
	\]	
	which is precisely~\eqref{eq:othercmkequiv} and, as a result, 
	proves~\eqref{eq:othercmk}.
\end{proof}\\

The formula~\eqref{eq:othercmk} is associated with another characterization of 
P$k$SSs, namely that provided by Definition~\ref{def:pkspanpkss}. In essence, 
the $k$-cosine measure of a family $\D$ is the minimum among all subsets of 
$\D$ with cardinality $|\D|-k+1$. A useful corollary of that result is that 
for any $k \ge 2$, one has
$$
	\cm{\D}=\cmk{\D}{1} \ge \cmk{\D}{2} \ge \cdots \ge \cmk{\D}{k}.
$$
In the next section, we will show how P$k$SSs built using OSPBs satisfy 
stronger properties associated with the $k$-cosine measure.

%%%%%%%%%%%%%%%%%%%%%%%%%%%%%%%%%%%%%%%%%%%%%%%%%%%%%%%%%%%%%%%%%%%%%%%%%%%%%%%
\subsection{Building positive \texorpdfstring{$k$}{TEXT}-spanning sets using OSPBs}
\label{subsec:pkssospb}
%%%%%%%%%%%%%%%%%%%%%%%%%%%%%%%%%%%%%%%%%%%%%%%%%%%%%%%%%%%%%%%%%%%%%%%%%%%%%%%

In this section, we describe how OSPBs can be used to generate positive 
$k$-spanning sets or positive $k$-bases with guarantees on their $k$-cosine 
measure. 

A first approach towards constructing positive $k$-spanning sets consists in 
duplicating the same PSS $k$ times, so that every direction remains even after 
taking out $k-1$ elements. However, such a strategy creates redundancy among 
the elements of the set. The purpose of this section is to introduce a more 
generic approach based on rotation matrices that allows for all vectors to be 
distinct.

\begin{proposition}
\label{prop:rotateospb}
	Let $\D_{n,s}$ be an OSPB of $\R^n$ with $s \in [\![1,n]\!]$. Let 
	$\Rm^{(1)},\dots,\Rm^{(k)}$ be $k$ rotation matrices in $\R^{n \times n}$, and 
	define $\D_{n,s}^{(j)}$ as the family of vectors obtained by applying 
	$\Rm^{(j)}$ to each vector in $\D_{n,s}$ for any $j \in [\![1,k]\!]$. Then, 
	the set $\D_{n,s}^{(1:k)}=\bigcup\limits_{j=1}^k \D_{n,s}^{(j)}$ is a positive 
	$k$-spanning set with
	\begin{equation}
	\label{eq:rotateospbcm}
		\cmk{\D_{n,s}^{(1:k)}}{k} \ge \cm{\D_{n,s}}.
	\end{equation}
\end{proposition}
\begin{proof}
	First, note that applying a rotation to every vector in a family does 
	not change its cosine measure since rotations preserve angles, hence 
	$\cm{\D_{n,s}^{(j)}}=\cm{\D_{n,s}}$ for every $j \in [\![1,k]\!]$.
	
	Consider a family $\setS \subset \D_{n,s}^{(1:k)}$ with 
	$|\setS| = \left|\D_{n,s}^{(1:k)}\right|-k+1$. By the pigeonhole principle, 
	this set must contain one of the $k$ positive bases obtained by rotation. 
	Letting $\D_{n,s}^{(j_\setS)}$ denote that positive basis, we have
	$\cm{\setS} \ge \cm{\D_{n,s}^{(j_\setS)}}=\cm{\D_{n,s}}$. 
	Using Theorem~\ref{thm:othercmk}, we obtain that
	\[
		\cmk{\D_{n,s}^{(1:k)}}{k}
		= \min_{\substack{\setS \subset \D_{n,s}^{(1:k)} \\
		|\set| = \left|\D_{n,s}^{(1:k)}\right|-k+1}} \cm{\setS} 
		\ge \cm{\D_{n,s}}.
	\]
	In particular, this proves $\cmk{\D_{n,s}^{(1:k)}}{k}>0$, hence this set 
	being positively $k$-spanning now follows from Theorem~\ref{thm:cmkpkss}.
\end{proof}\\

Several remarks are in order regarding Proposition~\ref{prop:rotateospb}. 
First, note that the result extends to any positive spanning set, but we focus 
on OSPBs since in that case~\eqref{eq:rotateospbcm} gives a lower bound on the 
$k$-cosine measure that can be computed in polynomial time. Moreover, when 
$\Rm_1=\dots=\Rm_k=\mI_n$, $i.e.$ when $k$ identical copies of $\D_{n,s}$ are 
used, the resulting family is a positive $k$-basis and 
relation~\eqref{eq:rotateospbcm} holds with equality. As a result, its 
$k$-cosine measure can be computed in polynomial time using 
Algorithm~\ref{algo:cmOSPB}. In general, however, the family generated through 
Proposition~\ref{prop:rotateospb} is not necessarily a positive $k$-basis. For 
instance, if $n=2$, setting $\D_{2}=\{\setI_2,-\one\}$, $k=2$, $\Rm_1=\mI_2$ 
and $\Rm_2 = \bbm 1/2 &-\sqrt{3}/2 \\ \sqrt{3}/2 &1/2 \ebm$ yields a family 
that is a positive $2$-spanning set but not a positive $2$-basis.

We now present a strategy based on rotations tailored to OSPBs. Recall that 
OSPBs can be decomposed into minimal positive bases over orthogonal subspaces. 
Our proposal consists in applying separate rotations to each of those 
minimal positive bases. Our key result is described in 
Theorem~\ref{thm:rotateminpb}, and shows that one can define a strategy for 
rotating a minimal positive basis to obtain a positive $k$-basis. The proof 
relies on two technical results, the first one being an intrinsic property of 
minimal positive bases.
\begin{lemma}\label{lem min pb strict sep}
	Let $\D_{\LL}=\{\vd_1,\dots,\vd_{\ell+1}\}$ be a minimal positive basis of 
	an $\ell$-dimensional subspace $\LL$ in $\R^n$. Then, there exists 
	$\vv_1,\dots,\vv_{\ell+1} \in \LL$ satisfying 
	\begin{equation}\label{eq:strict sep}
		\forall i \in [\![1, \ell+1]\!],\quad \vd_i^\top\vv_i > 0 
		\quad \text{and} \quad 
		\forall 1 \le j< i \le \ell+1, \quad \vd_i^\top\vv_{j}<0.
	\end{equation}
\end{lemma}
\begin{proof}
	By Lemma~\ref{lem:proppss}$\eqref{lem: PSS eq pos dotprod}$, the zero vector 
	$\zeron$ can be obtained by a positive linear combination of all elements 
	in $\D_{\LL}$. Without loss of generality, we rescale the elements of 
	$\D_\LL$ to ensure that $\sum\limits_{i=1}^{\ell+1}\vd_i=\zero_n$. Now, 
	since $\D_\LL$ is a minimal positive basis, the set 
	$\D_\LL\backslash\{\vd_{\ell+1}\}$ is a linear basis of $\LL$. Suppose 
	that we augment this set with $n-\ell$ vectors to form a basis $\B_n$ 
	of $\R^n$, and let $\Bm_n$ be a matrix representation of that basis 
	such that the first $\ell$ columns correspond to $\vd_1,\dots,\vd_{\ell}$. 
	Then, we have $\vd_i = \Bm_n \ve_i$ for any $i \in [\![1,\ell]\!]$, 
	while $\vd_{\ell+1} = -\Bm_n \left(\sum_{j=1}^{\ell} \ve_j\right)$, where we 
	recall that $\ve_1,\dots,\ve_{\ell}$ are the first $\ell$ vectors of 
	the canonical basis. 

	We now define the vectors 
	\[
		\left\{
			\begin{array}{lll}
				\vv_i &= &\Bm_n^{-\T} \left(-\sum_{j=1}^{\ell} 
				\ve_{j} + (\ell+1) \ve_i \right) \quad i \in [\![1,\ell]\!] \\
				\vv_{\ell+1} &= &-\Bm_n^{-\T} \sum_{j=1}^{\ell} \ve_j.
			\end{array}
		\right.
	\]	
	Using orthogonality of the coordinate vectors in $\R^n$, we have 
	that $\vd_i^\T \vv_i = \ell \|\ve_i\|^2 = \ell>0$ for 
	every $i \in [\![1,\ell]\!]$ and 
	$\vv_{\ell+1}^\T \vd_{\ell+1} = \left\|\sum_{j=1}^{\ell} \ve_j\right\|^2 >0$. 
	As a result, the first part of~\eqref{eq:strict sep} holds.
	
	In addition, for any $1 \le j < i \le \ell+1$, we obtain that 
	$\vd_i^\T \vv_{j} = -\|\ve_i\|^2=-1 < 0$ if $i < \ell+1$ and
	\[
		\vd_{\ell+1}^\T \vv_{j} 
		= \left\|\sum_{i=1}^{\ell} \ve_i\right\|^2 
		- (\ell+1) \|\ve_{j}\|^2
		= \sum_{i=1}^{\ell} \|\ve_i\|^2 - (\ell+1) \|\ve_{j}\|^2 
		= \ell - (\ell+1)=-1<0,
	\]
	thus the second part of~\eqref{eq:strict sep} also holds.
\end{proof}\\

Note that Lemma~\ref{lem min pb strict sep} is not equivalent to 
Definition~\ref{def:pi} as the inequalities in (\ref{eq:strict sep}) are 
strict.

Our second technical result uses the result of 
Lemma~\ref{lem min pb strict sep} to define a quantity characteristic of 
the angles between vectors in a minimal positive basis.

\begin{lemma}\label{lem:barrho}
	Suppose that $n\geq \ell\geq 2$ and let 
	$\D_{\LL}=\{\vd_1,\dots,\vd_{\ell+1}\}$ be a minimal positive basis of 
	a $\ell$-dimensional subspace $\LL$ in $\R^n$. Let $\vv_1,\dots,\vv_{\ell+1}$ 
	be vectors in $\LL$ that satisfy~\eqref{eq:strict sep}. For any pair 
	$(i,i') \in [\![1,\ell+1]\!]^2$, let $\vd_{i,i'}$ be the orthogonal 
	projection of $\vd_i$ on $\LL\cap\{\vv_{i'}\}^\perp$. Then, the 
	quantity 
	$\rho_{\D_{\LL}}:=\max\limits_{i,i' \leq \ell+1}
	\frac{\vd_i^\top\vd_{i,i'}}{\|\vd_i\|\|\vd_{i,i'}\|}$ lies in the 
	interval $[0,1)$.
\end{lemma}
\begin{proof}
	Since $\vd_{i,i'}$ is a projection of $\vd_i$ for any pair $(i,i')$, 
	we have that $\vd_i^\T \vd_{i,i'} \ge 0$, hence $\rho_{\D_{\LL}} \ge 0$. 
	Moreover, for any pair $1 \le i \le i' \le \ell+1$, we have 
	$\vd_i^\top \vv_{i'} \neq 0$ by~\eqref{eq:strict sep}, hence 
	$\vd_i \notin \LL \cap\{\vv_j\}^\perp$ and 
	$\tfrac{\vd_i^\top\vd_{i,i'}}{\|\vd_i\|\|\vd_{i,i'}\|} < 1$. As a 
	result, $\rho_{\D_{\LL}}< 1$, completing the proof.
\end{proof}\\

The quantity defined in Lemma~\ref{lem:barrho} is instrumental in defining 
suitable rotations matrices to produce a positive $k$-basis. The construction 
is described in Theorem~\ref{thm:rotateminpb}.

\begin{theorem}\label{thm:rotateminpb}
	Suppose that $n \ge \ell \ge 2$ and let 
	$\D_{\LL}=\{\vd_1,\dots,\vd_{\ell+1}\}$ be a minimal positive basis of an 
	$\ell$-dimensional subspace $\LL$ in $\R^n$. Let $\rho_{\D_{\LL}}$ be the 
	quantity defined in Lemma \ref{lem:barrho} for some vectors 
	$\vv_1,\dots,\vv_{\ell+1}$ satisfying~\eqref{eq:strict sep}. Finally, let 
	$\Rm^{(1)}_{\LL},\dots,\Rm^{(k)}_{\LL}$ be $k$ rotation matrices 
	in $\R^{n \times n}$ such that 
 	\begin{enumerate}[(i)]
 		\item\label{stayinL} $\Rm_\LL^{(j)}\vu=\vu, \quad$ for any 
 		$j \in [\![1,k]\!]$ and any vector $\vu \in \LL^\perp$,
 		\item\label{noredundancy} for any pair $(j,j') \in [\![1,k]\!]^2$ and 
 		any pair $(i,i') \in [\![1,\ell+1]\!]^2$,
		\[
			\frac{\left[ \Rm^{(j)}_{\LL}\vd_i\right]^\T 
			\Rm^{(j')}_{\LL}\vd_{i'}}
			{\|\Rm^{(j)}_{\LL}\vd_i\|\|\Rm^{(j')}_{\LL}\vd_{i'}\|} 
			= 1 
			\quad \Leftrightarrow \quad 
			j=j'\ \mbox{and}\ i=i',
		\]
 	\item\label{eq:smallangle}$\frac{\left[\Rm_{\LL}^{(j)}\vd_i\right]^\top\vd_i}
 	{\|\Rm_{\LL}^{(j)}\vd_i\|\|\vd_i\|}>\rho_{\D_{\LL}} 
 	\quad \text{for any pair} \quad (i,j) \in [\![1,\ell+1]\!]\times[\![1,k]\!]$.
 	\end{enumerate}
	Denote by $\D_{\LL}^{(j)}$ the family obtained by applying 
	$\Rm_{\LL}^{(j)}$ to each vector in $\D_{\LL}$. Then, the family
	$\dlk=\bigcup\limits_{j=1}^k \D_{\LL}^{(j)}$ is a positive $k$-basis 
	of $\LL$ with no identical elements.
\end{theorem}
\begin{proof}
	Owing to assumptions $(\ref{stayinL})$ and $(\ref{noredundancy})$, we know that 
	$\mathcal{D}_\LL^{(1:k)} $ is a subset of $\LL$ with pairwise distinct elements.
	Moreover, Proposition~\ref{prop:rotateospb} guarantees that $\dlk$ is a P$k$SS 
	for $\LL$. Therefore we only need to show the positive $k$-independence of this set.
	
	To this end, we consider an index $i \in [\![1,\ell+1]\!]$ and show that the 
	vector $\vv_i$ makes a positive dot product with exactly $k$ elements of 
	$\mathcal{D}_\LL^{(1:k)}$, these elements being 
	$\Rm_{\LL}^{(1)}\vd_i,\dots,\Rm_{\LL}^{(k)}\vd_i$. For any $j \in [\![1,k]\!]$, 
	we deduce from condition~$(\ref{eq:smallangle})$ in the theorem's statement 
	that 
	\[
		\frac{\left[\Rm_{\LL}^{(j)}\vd_i\right]^\top\vd_i}
		{\|\Rm_{\LL}^{(j)}\vd_i\|\|\vd_i\|}
		>\rho_{\D_{\LL}} \ge \frac{\vd_i^\top\vd_{i,j}}{\|\vd_i\|\|\vd_{i,j}\|},
	\]
	where $\vd_{i,j}$ is defined as in Lemma~\ref{lem:barrho}. Letting $\Theta(\vu,\vu')$ 
	denote the angle (with values in $[0,\pi]$) between two elements $\vu$ and $\vu'$ in 
	$\R^n$, the above property translates to 
	$\Theta(\Rm_{\LL}^{(j)}\vd_i,\vd_i)<\Theta(\vd_i,\vd_{i,j})$. Using now that 
	$\Theta(\vu,\vu') \le \Theta(\vu,\vu'') + \Theta(\vu'',\vu')$ for any vector 
	triplet $(\vu,\vu',\vu'') \in (\R^n)^3$, we obtain 
	\begin{equation*}
		\Theta(\Rm_{\LL}^{(j)}\vd_i,\vv_i) 
		\le \Theta(\Rm_{\LL}^{(j)}\vd_i,\vd_i)+\Theta(\vd_i,\vv_i) 
		< \Theta(\vd_{i,i},\vd_{i})+\Theta(\vd_i,\vv_i) 
		= \frac{\pi}{2},
	\end{equation*}
	where the last equality comes from the fact that $\vd_{i,i}$ and $\vv_i$ are 
	orthogonal vectors. 
	We have thus established that 
	$\vv_i^\top \Rm_{\LL}^{(j)}\vd_i>0$ for all $j \in [\![1,k]\!]$.
	
	It now remains to show that $\vv_i^\top\Rm_\LL^{(j)}\vd_{i'}\le 0$ for any 
	$i' \neq i$. First, using the same argument as above yields 
	$\Theta(\vv_{i},\vd_{i'})=\frac{\pi}{2}+\Theta(\vd_{i',i},\vd_{i'})$ as well as 
	$\Theta(\vv_{i},\Rm_{\LL}^{(j)}\vd_{i'})\ge 
	\Theta(\vv_{i},\vd_{i'})-\Theta(\vd_{i'},\Rm_{\LL}^{(j)}\vd_{i'})$. Applying 
	condition~$(\ref{eq:smallangle})$, we find that
	$\Theta(\vd_{i'},\Rm_{\LL}^{(j)}\vd_{i'})<\Theta(\vd_{i'},\vd_{i',i})$, which 
	finally leads to 
	\[
		\Theta(\vv_{i},\Rm_{\LL}^{(j)}\vd_{i'})
		\ge \Theta(\vv_{i},\vd_{i'})-\Theta(\vd_{i'},\Rm_{\LL}^{(j)}\vd_{i'})
		>\frac{\pi}{2},
	\]
	hence $\vv_i^\top\Rm_\LL^{(j)}\vd_{i'}\le 0$. Overall, we have established 
	that the vector $\vv_i$ makes a positive scalar product with exactly $k$ 
	vectors in $\dlk$. Since any vector $\vd \in \dlk$ is obtained by applying 
	a rotation to some element $\vd_i \in \dlk$, we conclude that $\dlk$ is 
	positively $k$-independent, proving the desired result.
\end{proof}\\

Note that Theorem~\ref{thm:rotateminpb} only applies when $\ell \ge 2$ (thus 
requiring $n \ge 2$ as well). Indeed, when the subspace $\LL$ has dimension $1$, 
none of the three conditions $(\ref{stayinL})$, $(\ref{noredundancy})$ and 
$(\ref{eq:smallangle})$ can be fulfilled.

\begin{remark}
\label{rk:condRotations}
	The first two conditions enforced on the rotation matrices in 
	Theorem~\ref{thm:rotateminpb} are designed to produce distinct vectors in 
	$\LL$ without affecting the orthogonal complement of $\LL$ in $\R^n$. Indeed, 
	Condition $(\ref{stayinL})$ guarantees that the rotation leaves any vector 
	orthogonal to $\LL$ invariant. Condition (ii) enforces that all vectors 
	produced by applying the rotations are distinct since the angle between two 
	different vectors has cosine less than 1. Finally, condition 
	$(\ref{eq:smallangle})$ ensures that the vectors are positively 
	$k$-independent. These conditions can be ensured through careful control of 
	the eigenvalues of those rotation matrices, according to the angles between 
	the vectors in $\D_{\LL}$.
\end{remark}

Theorem~\ref{thm:rotateminpb} provides a principled way of 
building positive $k$-bases from minimal positive bases. This approach naturally 
extends to OSPBs through their decomposition into orthogonal minimal positive 
bases. Using the technique described in Theorem~\ref{thm:rotateminpb} one can 
define rotation matrices that act on a single subspace from the decomposition of 
the OSPB without affecting the remaining subspaces thanks to orthogonality. The 
next corollary states the result.

\begin{corollary}
\label{cor:rotateospbdecomp}
	Let $\D_{n,s}$ be an OSPB and let $\D_{\LL_1},\dots,\D_{\LL_s}$ be a 
	decomposition~\eqref{eq:decompospb} of $\D_{n,s}$. Assume that for every $i\in[\![1,s]\!]$, $dim(\LL_i)>1$. Let $k\ge 1$ and 
	for every $i \in [\![1,s]\!]$, let $\Rm_{\LL_i}^{(1)},\dots,\Rm_{\LL_i}^{(k)}$ 
	be $k$ rotation matrices satisfying the properties (i), (ii) and (iii) from 
	Theorem~\ref{thm:rotateminpb} relative to $\D_{\LL_i}$. Then, the set
	$$
		\D_{n,s}^{(1:k)} = \bigcup_{j=1}^k \D_{n,s}^{(j)},
		\qquad \mbox{with} \quad 
		\D_{n,s}^{(j)} = \bigcup_{i=1}^s \D_{\LL_i}^{(j)}
		\quad \forall j \in [\![1,k]\!],
	$$
	is a positive $k$-basis with $\cmk{\D_{n,s}^{(1:k)}}{k} \ge \cm{\D_{n,s}}$. 
\end{corollary}

Corollary~\ref{cor:rotateospbdecomp} thus provides a useful strategy to design 
positive $k$-bases with guaranteed $k$-cosine measure, since a lower bound on 
that quantity is given by $\cm{\D_{n,s}}$ and that measure can be efficiently 
computed through Algorithm~\ref{algo:cmOSPB}. In particular, applying this 
strategy to a minimal positive basis yields a simple way of generating a 
positive $k$-basis with $k(n+1)$ vectors, as shown by 
Theorem~\ref{thm:rotateminpb} above. 

Similarly to the result of Theorem~\ref{thm:rotateminpb}, the result of 
Corollary~\ref{cor:rotateospbdecomp} does not apply to all OSPBs. In 
particular, the maximal OSPB $\{\In,-\In\}$ decomposes into $n$ minimal 
positive bases over one-dimensional subspaces, which precludes the application 
of Theorem~\ref{thm:rotateminpb}. Note however that 
Proposition~\ref{prop:rotateospb} still provides a way to compute P$k$SSs 
with cosine measure guarantees for such a maximal OSPB.

%%%%%%%%%%%%%%%%%%%%%%%%%%%%%%%%%%%%%%%%%%%%%%%%%%%%%%%%%%%%%%%%%%%%%%%%%%%%%%%
\section{Conclusion}
\label{sec:conc}
%%%%%%%%%%%%%%%%%%%%%%%%%%%%%%%%%%%%%%%%%%%%%%%%%%%%%%%%%%%%%%%%%%%%%%%%%%%%%%%

In this paper, we have studied two classes of positive spanning sets. On the 
one hand, we have investigated orthogonally structured positive bases, that 
possess a favorable structure in that they decompose into minimal positive 
bases over orthogonal subspaces. By exploiting this property, we described an 
algorithm that computes this structure as well as the value of the cosine 
measure in polynomial time, thereby improving an existing procedure for 
generic positive spanning sets. On the other hand, we have conducted a 
detailed study of positive $k$-spanning sets, a relatively understudied class 
of PSSs with resilient properties. We have provided several characterizations 
of these sets, including the generalization of the cosine measure through the 
$k$-cosine measure. We have also leveraged OSPBs to build positive 
$k$-spanning sets with guarantees on their $k$-cosine measure based on 
rotations.

Our results open the way for several promising areas of research. We believe 
that the cosine measure calculation technique can be further improved for 
positive bases by leveraging their decomposition, although the presence of 
critical vectors poses a number of challenges to overcome for dealing with 
positive bases that are not orthogonally structured. Such results would also 
improve the practical calculation of $k$-cosine measures. Finally, designing 
derivative-free optimization techniques based on positive $k$-spanning sets 
with guarantees on their $k$-cosine measure is a natural application of our 
study, and is the subject of ongoing work.

%\paragraph{Acknowledgments} We thank the Editor-in-Chief for handling our 
%paper, as well as an anonymous referee for their perceptive comments.

%%%%%%%%%%%%%%%%%%%%%%%%%%%%%%%%%%%%%%%%%%%%%%%%%%%%%%%%%%%%%%%%%%%%%%%%%%%%%%%
%%%%%%%%%%%%%%%%%%%%%%%%%%%%%%%%%%%%%%%%%%%%%%%%%%%%%%%%%%%%%%%%%%%%%%%%%%%%%%%
\bibliographystyle{siam}
\bibliography{bibliography}

\def\cprime{$'$}
\begin{thebibliography}{10}

\bibitem{Audet2011short}
{\sc C.~Audet}, {\em A short proof on the cardinality of maximal positive
  bases}, Optim. Letters, 5 (2011), pp.~191--194.

\bibitem{AudetHare2018derivative}
{\sc C.~Audet and W.~Hare}, {\em Derivative-Free and Blackbox Optimization},
  Springer Series in Operations Research and Financial Engineering, Springer
  International Publishing, 2018.

\bibitem{Conn2009}
{\sc A.~Conn, K.~Scheinberg, and L.~Vicente}, {\em Introduction to
  derivative-free optimization}, vol.~8, Siam, 2009.

\bibitem{Davis1954positive}
{\sc C.~Davis}, {\em Theory of positive linear dependence}, American Journal of
  Mathematics, 76 (1954), pp.~733--746.

\bibitem{Hare2018}
{\sc W.~Hare and G.~Jarry-Bolduc}, {\em Calculus identities for generalized
  simplex gradients: Rules and applications}, Submitted to SIAM Journal on
  Optimization,  (2018).

\bibitem{hare2020cm}
\leavevmode\vrule height 2pt depth -1.6pt width 23pt, {\em A deterministic
  algorithm to compute the cosine measure of a finite positive spanning set},
  Optim. Lett., 14 (2020), pp.~1305--1316.

\bibitem{hare2023nicely}
{\sc W.~Hare, G.~Jarry-Bolduc, and C.~Planiden}, {\em Nicely structured
  positive bases with maximal cosine measure}, Optim. Lett.,  (2023).

\bibitem{Hopcroft1973algograph}
{\sc J.~Hopcroft and R.~Tarjan}, {\em Algorithm 447: Efficient algorithms for
  graph manipulation}, Commun. ACM, 16 (1973), p.~372–378.

\bibitem{jarrybolducthesis}
{\sc G.~Jarry-Bolduc}, {\em Structures in Derivative Free Optimization:
  Approximating Gradients/ Hessians, and Positive Bases}, PhD thesis,
  University of British Columbia, 2023.
\newblock \url{https://circle.ubc.ca/}.

\bibitem{kolda2003directsearch}
{\sc T.~G. {Kolda}, R.~M. {Lewis}, and V.~Torczon}, {\em Optimization by direct
  search: {N}ew perspectives on some classical and modern methods}, SIAM Rev.,
  45 (2003), pp.~385--482.

\bibitem{lewis1996rank}
{\sc R.~M. Lewis and V.~Torczon}, {\em Rank ordering and positive bases in
  pattern search algorithms.}, tech. rep., Institute for computer applications
  in science and engineering, Hampton VA, 1996.

\bibitem{Marcus1981}
{\sc D.~A. {Marcus}}, {\em Minimal positive 2-spanning sets of vectors},
  Proceedings of the AMS, 82 (1981), pp.~165--172.

\bibitem{Marcus1984}
{\sc D.~A. Marcus}, {\em Gale diagrams of convex polytopes and positive
  spanning sets of vectors}, Discrete Applied Mathematics, 9 (1984),
  pp.~47--67.

\bibitem{marsden2013eigenvalues}
{\sc A.~Marsden}, {\em Eigenvalues of the {L}aplacian and their relationship to
  the connectedness of a graph}, University of Chicago, REU,  (2013).

\bibitem{mckinney1962}
{\sc R.~McKinney}, {\em Positive bases for linear spaces}, Transactions of the
  American Mathematical Society, 103 (1962), pp.~131--148.

\bibitem{naevdal2019}
{\sc G.~Naevdal}, {\em Positive bases with maximal cosine measure},
  Optimization Letters, 13 (2019), pp.~1381--1388.

\bibitem{Regis2015}
{\sc R.~Regis}, {\em The calculus of simplex gradients}, Optimization Letters,
  9 (2015), pp.~845--865.

\bibitem{Regis2016}
{\sc R.~Regis}, {\em On the properties of positive spanning sets and positive
  bases}, Optimization and Engineering, 17 (2016), pp.~229--262.

\bibitem{Regis2021}
{\sc R.~Regis}, {\em On the properties of the cosine measure and the uniform
  angle subspace}, Comput. Optim. Appl., 78 (2021), pp.~915--952.

\bibitem{roberts2022subspace}
{\sc L.~Roberts and C.~W. {Royer}}, {\em Direct search based on probabilistic
  descent in reduced spaces}.
\newblock arXiv:2204.01275v2, 2023.

\bibitem{Romanowicz1987}
{\sc Z.~Romanowicz}, {\em Geometric structure of positive bases in linear
  spaces}, Applicationes Mathematicae, 19 (1987), pp.~557--567.

\bibitem{Wotzlaw2009}
{\sc R.~Wotzlaw}, {\em Incidence graphs and unneighborly polytopes}, PhD
  thesis, Technischen Universit\"at Berlin, 2009.

\end{thebibliography}
%%%%%%%%%%%%%%%%%%%%%%%%%%%%%%%%%%%%%%%%%%%%%%%%%%%%%%%%%%%%%%%%%%%%%%%%%%%%%%%
%%%%%%%%%%%%%%%%%%%%%%%%%%%%%%%%%%%%%%%%%%%%%%%%%%%%%%%%%%%%%%%%%%%%%%%%%%%%%%%

\end{document}